\def\Date{August 16, 2009}

\magnification=\magstephalf
\advance\vsize14mm
\frenchspacing\parskip4pt plus 1pt

\newread\aux\immediate\openin\aux=\jobname.aux 
\ifeof\aux\message{ <<< Run TeX a second time >>> }
\batchmode\else\input\jobname.aux\fi\closein\aux
\input \jobname.def\immediate\openout\aux=\jobname.aux

\draftfalse\footline={\hss}  

\headline={\ifnum\pageno>1\eightrm\ifodd\pageno\hfill~~ Fels and Kaup,
CR-tube realizations
\hfill{\tenbf\folio}\else{\tenbf\folio}\hfill Fels and Kaup,
maximal abelian subalgebras~~ \hfill\fi\else\hss\fi}

\centerline{\Gross Local tube realizations of CR-manifolds and maximal}
\medskip
\centerline{\Gross abelian subalgebras}
\medskip\bigskip
\centerline{{GREGOR FELS~~~and~~~WILHELM KAUP}}

{\bigskip\bigskip\openup-2pt\narrower\noindent\ninebf
Abstract\ninerm~~ For every real-analytic CR-manifold $M$ we give
necessary and sufficient conditions that $M$ can be realized in a
suitable neighbourhood of a given point $a\in M$ as a tube submanifold
of some $\CC^{r}$.  We clarify the question of the `right' equivalence
between two local tube realizations of the CR-manifold germ $(M,a)$ by
introducing two different notions of affine equivalence. One of our
key results is a procedure that reduces the classification of
equivalence classes to a purely algebraic manipulation in terms of Lie
theory.\par\noindent\ninebf Mathematics Subject
Classification\ninerm~~ Primary 32V05; Secondary 32V40, 32M25,
17B66\par}

\bigskip\medskip

\KAP{Introduction}{Introduction}

Among all CR-submanifolds of $\CC^{r}$ a special class is formed by
the {\sl tube submanifolds,} that is, by real submanifolds of the form
$$T_{F}=\RR^{r}+iF\Leqno{GH}$$ with $F$ an arbitrary submanifold of
$\RR^{r}$, called the {\sl base} of $T_{F}$.  CR-manifolds of this
type play a fundamental role in CR-geometry as they often serve as
test objects. In addition, the interplay between real geometric
properties of the base $F$ and CR-properties of the associated tube
$T_{F}$ are quite fruitful.  An early example of this interplay is
well known in the case of open tube submanifolds: {\sl The tube domain
$T_{F}\subset\CC^{r}$ is holomorphically convex if and only if the
(open) base $F\subset\RR^{r}$ is convex in the elementary sense.}
Clearly, in the context of CR-geometry, domains in $\CC^{r}$ are not
of interest. In fact, we will mainly consider CR-manifolds $M=(M,\H
M,J)$ which are holomorphically nondegenerate, i.e., $\xi=0$ is the
only holomorphic vector field on $M$, which is a section in the
subbundle $\H M$.  We note in passing that in the tube situation the
general case can be reduced to the nondegenerate one as every such
CR-manifold is locally a direct product of some $\CC^k$ and a
holomorphically nondegenerate CR-manifold.

For instance, interesting examples of holomorphically
nondegenerate  tube submanifolds are obtained  as
follows: Let $\Omega\subset\RR^{r}$ be an open convex cone such that
the corresponding tube domain $T_{\Omega}\subset\CC^{r}$ is
biholomorphically equivalent to an irreducible bounded symmetric
domain. Then the group
$G=\GL(\Omega):=\{g\in\GL(r,\RR):g(\Omega)=\Omega\}$ acts transitively
on $\Omega$ and for every non-open $G$-orbit $F\subset\RR^{r}$ with
$F\ne\{0\}$ the corresponding tube $T_{F}$ is Levi degenerate but
still is holomorphically nondegenerate \Lit{KAZT}. The example of
lowest possible dimension occurs with the future cone
$\Omega=\{x\in\RR^{3}:x_{3}>\sqrt{x_{1}^{2}+x_{2}^{2}}\}$ in
3-dimensional space-time and
$F=\{x\in\RR^{3}:x_{3}=\sqrt{x_{1}^{2}+x_{2}^{2}}>0\}$ the future
light cone. The future light cone tube $T_{F}$ has been studied by
many authors and has remarkable properties, compare \Lit{FEKA} and the
references therein. Until recently, this tube manifold $T_{F}$ was, up
to local CR-isomorphy, the only known example of a 5-dimensional Levi
degenerate, holomorphically nondegenerate and locally homogeneous
CR-manifold. A full classification of CR-manifolds of this type could
be obtained in \Lit{FKAU} -- surprisingly all possible examples turned
out to be locally representable  as tube manifolds.

Since tube manifolds are quite easy to deal with it is of interest to
decide whether a given CR-manifold $M$ is CR-isomorphic, at least
locally around a given point $a\in M$, to a tube submanifold of some
$\CC^{r}$. Another question is how many `different' tube realizations
a given CR-manifold germ does admit.  In the particular case of
spherical hypersurfaces the following result has been obtained in
\Lit{DAYA} by solving a certain partial differential equation coming
from the Chern-Moser theory \Lit{CHMO}: {\sl For every $r\ge2$ there
exist, up to affine equivalence, precisely $r+2$ closed smooth tube
submanifolds of $\CC^{r}$ that are locally CR-isomorphic to the
euclidian sphere $S^{2r-1}\subset\CC^{r}$.} In \Lit{ISMI}, \Lit{ISAE}
the same method has been used for a certain more general class of
CR-flat manifolds. All the above results rely on Chern-Moser theory
and therefore only apply to CR-manifolds that are Levi nondegenerate
and of hypersurface type.

In this note we use a different method that applies to all
CR-manifolds (for simplicity we work in the category of real-analytic
CR-manifolds). This method is more algebraic in nature and starts from
the following simple observation: A real submanifold $M\subset\CC^{r}$
is tube \Ruf{GH} if and only if $M$ is invariant under all real
translations $z\mapsto z+v$ with $v\in\RR^{r}$.  In particular,
$\7g:=\hol(M,a)$, the Lie algebra of all (germs of real-analytic)
infinitesimal CR-transformations at $a$, contains the abelian Lie
subalgebra induced by the above translations.  Therefore it is not
unexpected that every tube realization of an arbitrarily given
CR-manifold germ $(M,a)$ is strongly related to a certain abelian Lie
subalgebra $\7v$ of $\hol(M,a),$ see Prop. \ruf{UR} and Prop. \ruf{KR}
for precise statements.

In a slightly different form the Lie algebra $\7v$ has already been
used in \Lit{BART} for the characterization of tube manifolds (in fact
more generally in the context of abstract smooth CR-manifolds and the
solution of the local integrability problem for rigid CR-manifolds; on
the other hand we do not need to assume that the evaluation map
$\varepsilon_a:\7v\to T_aM$ is injective).  But, in contrast to
\Lit{BART} our intentions are completely different - we mainly focus
on the question how may `essentially' different tube realizations of a
given CR-manifold germ $(M,a)$ do exist.  This question of equivalence
for different local tube realizations of a given CR-manifold is a bit
more subtle than it might appear at the first glance.  We introduce
two different notions of equivalence to which we refer accordingly as
to the `strict' and the `coarse' affine equivalence. Our impression is
that the latter one is more appropriate in the context of local tube
realizations.

In Section \ruf{realizations} we give necessary and sufficient
conditions for an abelian subalgebra $\7v\subset\7g$ to give a local
tube realization of $(M,a)$. This characterization also includes for
every $\7v$ an easy to compute canonical form of a local
CR-isomorphism to the corresponding tube realization of $(M,a)$.  It
is also shown that any two local tube realizations of the germ $(M,a)$
are affinely equivalent (in the strict sense) if and only if the
corresponding abelian subalgebras $\7v,\7v'\subset\7g$ are conjugate
with respect to the stability group $\Aut(M,a)$.

The `coarse' equivalence relation for tube realizations of the germ
$(M,a)$ is, roughly speaking, defined as follows: Two tube
realizations $(T,c),(T',c')$ of $(M,a)$ in $\CC^{r}$ are considered to
be equivalent in this broader sense if the representing tube
submanifolds $T,T'\subset\CC^{r}$ can be chosen in such a way that
$T'=g(T)$ for some affine isomorphism $g$ on $\CC^{r}$ (that is, {\sl
without} requiring $c'=g(c)$ in addition). 

While it is not surprising that the existence of a tube realization
for $(M,a)$ is closely related to the existence of a certain `big'
abelian Lie subalgebra of $\7g=\hol(M,a)$, it is not at all clear what
the relation between various tube realizations and the corresponding
abelian subalgebras in $\hol(M,a)$ should be. One of our main results
is then obtained in Section \ruf{The}, where we introduce the subgroup
$\Glob(M,a)\subset\Aut(\7g)$ and show for a large class of
CR-manifolds $M$ that the local tube realizations of $(M,a)$ are
equivalent in the coarser sense if and only if the corresponding
abelian subalgebras $\7v,\7v'$ are conjugate with respect to the group
$\Glob(M,a)$.

In Sections \ruf{sphere} and \ruf{Further} we apply our general theory
to some concrete cases.  For instance, we relate the results from
\Lit{DAYA} with our algebraic point of view, and identify the various
abelian subalgebras of $\hol(S^{2r-1},a)$, $S^{2r-1}\subset\CC^{r}$
the standard sphere, which correspond to various
defining equations in \Lit{DAYA}.

In the last two sections we generalize the notion of a tube
submanifold to the notion of a Siegel submanifold.  This is motivated
by the well known fact that every bounded homogeneous domain can be
realized as a Siegel domain, thus giving a lot of additional insight
to the structure of those domains.  In the forthcoming paper
\Lit{FLKA} our method will be applied to the class of all Levi
non-degenerate real hyperquadrics in $\CC^{r}$ in order to obtain a
full algebraic characterization of local tube realizations in such
cases.

\vfil\eject

\KAP{Preliminaries}{Preliminaries and notation}

{\bf\noindent Abstract CR-manifolds.} A triple $(M,\H M,J)$ is called
an ({\sl abstract}) {\sl CR-manifold} (CR stands for Cauchy-Riemann)
if $M$ is a (connected if not stated otherwise explicitly) smooth
manifold, $\H M$ is a smooth subbundle of its tangent bundle $\T M$
and $J$ is a smooth bundle endomorphism of $\H M$ with
$J^{2}=-\id$. For simplicity we often write just $M$ instead of $(M,\H
M,J)$. For every $a\in M$ the restriction of $J$ to the linear
subspace $\H_{a}M\subset \T_{a}M$ makes $\H_{a}M$ to a complex vector
space, we call it the {\sl holomorphic tangent space to $M$ at $a$}
(in the literature $\H_{a}M$ is also called the {\sl complex
tangent space} and denoted by $\T^{\,c}_{a}M$).  Its complex dimension
is called the {\sl CR-dimension} and the real dimension of
$\T_{a}M/\H_{a}M$ is called the {\sl CR-codimension} of $M$. With
$M=(M,\H M,J)$ also $M\conj:=(M,\H M,-J)$ is a CR-manifold; we call it
the {\sl conjugate} of $M$.

A smooth map $g:M\to M'$ between two CR-manifolds
is called CR if for every
$a\in M$ and $a':=g(a)$ the differential $dg_{a}:\T_{a}M\to \T_{a'}M'$
maps the corresponding holomorphic subspaces in a complex linear way
to each other. Also, $g$ is called {\sl anti-CR} if $g$ is CR
as a map $M\conj\to M'$.

For every smooth vector field $\xi$ on $M$ and every $a\in M$ we
denote by $\xi_{a}\in \T_{a}M$ the corresponding tangent vector at
$a$. Furthermore, $\xi$ is called an {\sl infinitesimal
CR-transformation} of $M$ if the corresponding local flow on $M$
consists of CR-transformations. With $\xi,\eta$ also the usual bracket
$[\xi,\eta]$ is an infinitesimal CR-transformation.

It is obvious that every smooth manifold $M$ can be considered as a
CR-manifold with CR-dimension 0 (these are called the {\sl totally real}
CR-manifolds). The other extreme is formed by the CR-manifolds with
CR-codimension 0, these are precisely the almost complex
manifolds. Among the latter the integrable ones play a special role,
the complex manifolds. CR-mappings between complex manifolds are
precisely the holomorphic mappings.

{\bf\noindent CR-manifolds in this paper} are understood to be those
 $M=(M,\H M,J)$ that are real-analytic and integrable in the following
 sense: $M$ is a real-analytic manifold and there is a complex
 manifold $Z$ such that $M$ can be realized as a real-analytic
 submanifold $M\subset Z$ with $\H_{a}M=\T_{a}M\cap\,i\T_{a}M$ and
 $J(\xi)=i\xi$ for every $a\in M$, $\xi\in \H_{a}M$, where $\T_{a}M$
 is considered in the canonical way as an $\RR$-linear subspace of the
 complex vector space $\T_{a}Z$.  This notion of integrability is
 equivalent to the vanishing of the restricted Nijenhuis tensor. We
 refer to \Lit{BOGG} or \Lit{BERO} for further details.  The embedding
 $M\subset Z$ above can always be chosen to be {\sl generic}, that is,
 $\T_{a}Z=\T_{a}M+\,i\T_{a}M$ for all $a\in M$. In that case the
 (connected) complex manifold $Z$ has complex dimension (CR-$\dim
 M+\,$CR-$\codim M$).

CR-isomorphisms between CR-manifolds are always understood to be
analytic in both directions. In particular, $\Aut(M)$ is the group of
all (bianalytic) CR-automorphisms of $M$ and
$\Aut_{a}(M):=\{g\in\Aut(M):g(a)=a\}$ is the isotropy subgroup at the
point $a\in M$. With $\Aut(M,a)$, also called th {\sl stability group
at $a$}, we denote the group of all CR-automorphisms of the manifold
germ $(M,a)$. Then $\Aut_{a}(M)$ can be considered in a canonical way
as a subgroup of $\Aut(M,a)$.

With $\hol(M)$ we denote the space of all real-analytic infinitesimal
transformations of the CR-manifold $M$ and with $\hol(M,a)$ the space
of all germs at $a\in M$ of vector fields $\xi\in\hol(N)$ where $N$
runs through all open connected neighbourhoods of $a$ in $M$. Then
$\hol(M)$ as well as every $\hol(M,a)$ together with the bracket
$[\;,\;]$ is a real Lie algebra (of possibly infinite dimension).  The
canonical restriction mapping $\rho_{a}:\hol(M)\to\hol(M,a)$ is an
injective homomorphism of Lie algebras. Every isomorphism
$g:(M,a)\to(M',a')$ of CR-manifold germs induces in a canonical way a
Lie algebra homomorphism $g_{*}:\hol(M,a)\to\hol(M',a')$. Its inverse
is the pull back $g^{*}$. Clearly, $g\mapsto g_{*}$ defines a group
homomorphism $\Ad:\Aut(M,a)\to\Aut(\hol(M,a))$.

A vector field $\xi\in\hol(M)$ is called {\sl complete} on $M$ if the
corresponding local flow extends to a one-parameter group
$\RR\to\Aut(M)$. The image of $1\in\RR$ is denoted by $\exp(\xi)$. In
this sense we have the exponential map $\exp:\aut(M)\to\Aut(M)$, where
$\aut(M)$ is the set of all {\sl complete} $\xi\in\hol(M)$. In
general, $\aut(M)\subset\hol(M)$ is neither a linear subspace nor
closed under taking brackets. But, if there exists a Lie subalgebra
$\7g\subset\hol(M)$ of finite dimension with $\aut(M)\subset\7g$, then
$\aut(M)$ itself is a Lie subalgebra \Lit{PALA} and on $\Aut(M)$
there exists a unique Lie group structure (in general not connected)
such that $\exp$ is a local diffeomorphism in a neighbourhood of
$0\in\aut(M)$. Furthermore, the map $\Aut(M)\times M\to M$, $(g,a)\mapsto
g(a)$, is real-analytic.

In case $M$ is generically embedded as a real-analytic CR-submanifold
of a complex manifold $Z$ then a vector field $\xi$ on $M$ is in
$\hol(M)$ if and only if $\xi$ has an extension $\tilde\xi$ to a
holomorphic vector field on a suitable open neighbourhood $U$ of $M$
in $Z$ (that is, $\tilde\xi$ is a holomorphic section over $U$ in its
tangent bundle $\T U\,$). The Lie algebras $\hol(Z)$ and $\hol(Z,a)$
are complex Lie algebras and $\7g:=\hol(M,a)$ is in a canonical way a
real subalgebra of $\hol(Z,a)$. The CR-manifold germ $(M,a)$ is called
{\sl holomorphically nondegenerate} if $\7g$ is totally real in
$\hol(Z,a)$, that is, $\7g\cap\, i\!\7g=\{0\}$.  In this case there is
a unique antilinear Lie algebra automorphism $\sigma$ of
$\7g^{\CC}:=\7g+\,i\7g\subset\hol(Z,a)$ with
$\7g=\Fix(\sigma)$. Clearly, real Lie subalgebras of $\7g$ and
$\sigma$-invariant complex Lie subalgebras of $\7g^{\CC}$ are in a
natural 1-1-correspondence.

In general, a vector field $\xi\in\hol(M)$ only can be integrated to a
{\sl local} 1-parameter group of CR-transformations $g_{t}$
that we also denote by $\exp(t\xi)$. The reason for this notation in
the analytic case is the following: To every $a\in M$ and every open
neighbourhood $W$ of $a\in Z$ there is a further open neighbourhood
$U\subset W$ of $a\in Z$ and an $\epsilon>0$ such that the $g_{t}$ are
defined as holomorphic mappings $U\to W$ for $|t|<\epsilon$ and
satisfy for every holomorphic mapping $f:W\to\CC^{n}$ the formula
$$f\circ g_{t}|_{U}=\sum_{k=0}^{\infty}{1\over
k!}(t\xi)^{k}(f|_{U})\,.$$ In particular, if $f$ gives a local chart
for $Z$ around $a$ then the $g_{t}$ on $U$ can be recovered from the
right side of this formula.

\medskip

\Lemma{LU} Let $Z$ be a connected complex manifold of dimension $n$
and $\7e\subset\hol(Z,a)$ an abelian complex Lie subalgebra with
$\epsilon_{a}(\7e)=\T_{a}Z$, where $\epsilon_{a}$ is the evaluation map
$\xi\mapsto\xi_{a}$. Then $\epsilon_{a}$ induces a complex linear
isomorphism from $\7e$ onto $\T_{a}Z$. In particular, $\7e$ also has
dimension $n$ and is maximal abelian in $\hol(Z,a)$.

\Proof Let $\eta\in\7e$ be an arbitrary element with $\eta_{a}=0$. We
have to show $\eta=0$. Fix a linear subspace $\7a\subset\7e$ such that
$\epsilon_{a}:\7a\to\T_{a}Z$ is an isomorphism. We may assume that
$\7a\subset\hol(U)$ for some open neighbourhood $U\subset Z$ of $a$
and also that every $z\in U$ is of the form $z=\exp(\xi)(a)$ for some
$\xi\in\7a$. For every such $z$ then $[\eta,\xi]=0$ implies
$\exp(t\eta)(z)=\exp(\xi)\exp(t\eta)(a)=\exp(\xi)(a)=z$ for $|t|$
small, that is, $\eta=0$.\qed

For the sake of clarity we mention that in case $n=\dim Z\ge 2$
 there exist {\sl abelian} subalgebras $\7e\subset\hol(Z,a)$ of {\sl
 arbitrary} dimension. However, in general these do not span
 $\T_aZ$.

\medskip The CR-manifold $M$ is called {\sl homogeneous} if the group
$\Aut(M)$ acts transitively on $M$. Also, $M$ is called {\sl locally
homogeneous} if for every $a,b\in M$ the manifold germs $(M,a)$,
$(M,b)$ are CR-isomorphic. By \Lit{ZAIT} this is equivalent to
$\epsilon_{a}(\hol(M,a))=\T_{a}M$ for every $a\in M$. The
CR-manifold $M$ is called {\sl minimal} if every smooth submanifold
$N\subset M$ with $\H_{a}M\subset \T_{a}N$ for all $a\in N$ is already
open in $M$.

\medskip For later use (Proposition \ruf{ZB}) we state

\Lemma{BR} Let $Z$ be a complex manifold and $M\subset Z$ a (connected
real-analytic) generic and minimal CR-submani\-fold. Then for every 
closed complex-analytic subset $A\subset Z$ the set $M\backslash A$ is
connected.
 
\Proof We first show that the proof of the Lemma can be reduced to the
case when $A\subset Z$ is non-singular.  Indeed, there is an integer
$k\ge1$ and a descending chain $A=A_{0}\supset\cdots\supset
A_{k}=\emptyset$ of analytic subsets such that $A_{j}$ is the singular
locus of $A_{j-1}$ for all $j=1,\dots,k$. Put $M_{j}:=M\backslash
A_{j}$. Then $A_{j-1}\backslash A_{j}$ is analytic in
$Z_{j}:=Z\backslash A_{j}$ and
$M_{j-1}=M_{j}\backslash(A_{j-1}\backslash A_{j})$. Therefore it
suffices to show inductively that $M=M_{k},M_{k-1},\dots,M_{0}$ all
are connected. For the rest of the proof we therefore assume that $A$
is nonsingular and also, contrary to the claim of the Lemma, that
$M\backslash A$ is disconnected. Notice that this implies
$$\T_aM\cap\T_aA\;\ne\;\T_aM\Steil{for all}a\in M\cap A\,,\leqno{(*)}$$
since otherwise $M\backslash A=\emptyset$ would be connected as a
consequence of $\T_{a}Z=\T_{a}M+i\T_{a}M\subset\T_{a}A\subset\T_{a}Z$.

The intersection $S:=A\cap M$ is a real-analytic set. Again, there is
an integer $r\ge1$ and a descending chain $S=S_{0}\supset\cdots\supset
S_{r}=\emptyset$ of real-analytic subsets such that $S_{j}$ is the
singular locus of $S_{j-1}$ for all $j=1,\dots,r$. Choose $j\le r$
minimal with respect to the property that $M\backslash S_{j}$ is
connected. Then $j>0$ by the above assumption and $M\backslash
S_{j-1}=(M\backslash S_{j})\backslash (S_{j-1}\backslash S_{j})$ is
disconnected.  In particular, also $(M\backslash S_{j})\backslash N$
is disconnected, where we denote by $N$ the union of all connected
components of $(S_{j-1}\backslash S_{j})$ that have codimension 1 in
$M$. Clearly, $(*)$ improves to $$\T_aM\cap\T_aA\;=\;\T_aN\Steil{for
all}a\in N\,.\leqno{(**)}$$ Since $M$ is minimal by assumption there
exists an $a\in N$ with $\H_aM\not\subset \T_{\!a}N$ and hence with
$\H_aM\not\subset \T_{\!a}A$ by $(**)$. Since $\H_aM$ and $\T_aA$ are
complex linear subspaces, there is a linear subspace $V\subset
\H_aM\subset \T_{\!a}M$ of real dimension $\ge 2$ with $\H_aM=V\oplus
(\H_aM\cap\, \T_aA)$.  But then $V\cap\,\T_aN= V\cap(\T_aM\cap
\T_aA)=0$ gives a contradiction since $T_{a}N$ is a real hyperplane in
$T_{a}M$. This shows that $M\backslash A$ cannot be assumed to be
disconnected, and the proof is complete.\qed

Notice that the assumption on $M$ in Lemma \ruf{BR} is automatically
satisfied if $M$ is of hypersurface type and has nowhere vanishing
Levi form. Indeed, if $M$ is a hypersurface and is not minimal in
$a\in M$ then the Levi form of $M$ at $a$ vanishes.

{\bf\smallskip\noindent Convention for notating vector fields.}  
In this paper we do not
need the complexified tangent bundle $TM\otimes_{\RR}\CC$ of $M$. All
vector fields occurring here correspond to `real vector fields'
elsewhere. In particular, if $E$ is a complex vector space of finite
dimension and $U\subset E$ is an open subset then the vector fields
$\xi\in\hol(U)$ correspond to holomorphic mappings $f:U\to E$, and the
correspondence is given in terms of the canonical trivialization
$TU\cong U\times E$ by identifying the mapping $f$ with the vector
field $\xi=(\id_{U},f)$. To have a short notation we also write
$$\xi=f(z)\dd z{}\,.$$ As soon as the vector field $\xi=f(z)\dd z{}$
is considered as differential operator, special caution is necessary:
$\xi$ applied to the smooth function $h$ on $U$ is $\xi h=f(z)\dd
z{}h+\overline f(z)\dd{\overline z} h$. We therefore stress again that
we write
$$\xi=f(z)\dd z{}\Steil{ instead of }\xi=f(z)\dd z{}+\overline
f(z)\dd{\overline z}{}\hbox{~~elsewhere}\,,$$ and this convention will
be in effect allover the paper.

\KAP{manifolds}{Tube manifolds}

Throughout this section let $V$ be a real vector space of finite
dimension and $E:=V\oplus iV$ its complexification.  For every
(connected and locally closed) real-analytic submanifold $F\subset V$
the manifold
$$T:=T_{F}:=V+iF\;\subset\; E$$ is a CR-submanifold of $E$, called the
{\sl tube} over the {\sl base} $F$. Obviously, a real-analytic
submanifold $M\subset E$ is a tube in this sense if and only if
$M+V=M$. Tubes form a very special class of CR-manifolds. For
instance, $\Aut(T)$ contains the following abelian translation group
isomorphic to the vector group $V$
$$\Gamma:=\{z\mapsto z+v:v\in V\}\,.$$ Since $T=\Gamma(iF)$ it is
enough to study the local CR-structure of the tube $T$ only at points
$ia\in iF\subset T$. For these
$$\T_{ia}T=V\oplus i\T_{a}F\steil{and}\H_{ia}T=\T_{a}F\oplus
i\T_{a}F\subset E$$ is easily seen. In particular, $T$ is generic in
$E$. For every further tube $T'=V'+iF'$ in a complex vector space
$E'=V'\oplus iV'$ with $F'\subset V'$ every real affine mapping
$g:V\to V'$ with $g(F)\subset F'$ extends to a complex affine mapping
$\tilde g:E\to E'$ with $\tilde g(T)\subset T'$ and thus gives a
CR-map $T\to T'$. Therefore, $F$ (locally) being affinely homogeneous
implies that the tube $T$ is (locally) CR-homogeneous. The converse is
not true in general.

\Lemma{UO} Suppose that $T=V+iF$ is a tube submanifold of the complex
vector space $E=V\oplus iV$ and that $a\in T$ is an arbitrary
point. Then the following conditions are equivalent.  \0 $T$ is of
finite type at $a$. \1 $T$ is minimal at $a$. \1 The smallest
affine subspace of $V$ containing $F$ is $V$ itself.  \Proof It is
enough to show the implication \To31. We therefore assume (iii) and
identify $E=\CC^{n}$ with $\RR^{n}\times\RR^{n}$ in the standard way via
$(x+iy)\cong(x,y)$. Without loss of generality we assume that $T$
contains the origin of $E$ and is given in a suitable
neighbourhood of it by real-analytic equations
$$y_{j}=f_{j}(y_{1},\dots,y_{k}),\quad k<j\le n,$$ where every $f_{j}$
vanishes of order $\ge2$ at the origin of $\RR^k$. The assumption
(iii) implies that the germs of the functions $f_{k+1},\dots,f_{n}$ at
$0\in\RR^k$ are linearly independent. For all $1\le\ell,m\le k$ the
vector fields
$$\xi^\ell:=-\dd{x_\ell}+\sum_{j>k}\dD{f_{\!j}}{y_\ell}\dd{x_j}
\Steil{and}\eta^m:=\dd{y_m}-\sum_{j>k}\dD{f_{\!j}}{y_m}\dd{y_j}$$
(expressed in the real coordinates $(x,y)$ of $E$) are sections in the
holomorphic subbundle $\H T$ over the tube manifold $T$. Also, for
every multi-index $\nu=(\nu_1,\dots,\nu_k)\in\NN^k$ with
$|\nu|:=\nu_1+\dots+\nu_k\ge1$ and every $\ell=1,\dots,k$ we have
$$(\ad\eta^1)^{\nu_1}(\ad\eta^2)^{\nu_2}\cdots(\ad\eta^k)^{\nu_k}
\xi^\ell\;=\;\sum_{j>k}\Big(\dD{^{|\nu|}}{y^\nu}\big(\dD{f_{\!j}}{y_\ell}
\big)\Big)\dd{x_j}\,.\Leqno{YO}$$ Denote by $S\subset\T_{0}T$ the
linear subspace spanned by $\H_{0}T$ and all vector fields \Ruf{YO}.
Assume that
there exists a non-trivial linear form $\lambda$ on $\T_{0}T$ with
$\lambda(S)=0$ and put $f:=\sum_{j>k}f_{j}$ with
$d_{j}:=\lambda(\dd{z_{j}}\!)\,$. Then $d_{j}\ne0$ for some $j>k$,
that is, $f\not\equiv0$. On the other hand, \Ruf{YO} shows that all
partial derivatives of $f$ of order $\ge2$ vanish. By choice of the
functions $f_{j}$ also all partial derivatives of $f$ of order $<2$
vanish, a contradiction. Therefore $S=\T_{0}T$ and (i) must hold.\qed

\Proposition{} Suppose that $T=V+iF$ is a tube submanifold of the
complex vector space $E=V\oplus iV$ and suppose, without loss of
generality, that $T$ contains the origin of $E$. Then there exist
complex linear subspaces $E',E''$ of $E$ and tube submanifolds
$T'\subset E'$, $T''\subset E''$ with the following properties: \0
$T'$ is an $\RR$-linear subspace of $E'$ with $E'=T'+iT'$.\1
$T''$ is holomorphically nondegenerate and of finite type at every
point.\1 $E=E'\oplus E''$ and $T$ is open in $T'+T''$. \par\noindent

\Proof By Lemma \ruf{UO} we assume without loss of generality that $V$
is the linear span of $F$. We then verify the claim with
$T'=E'$.\nline For every $a\in T$ put $\7h_{a}:=\hol(T,a)\cap
i\hol(T,a)$ and $E_{a}:=\epsilon_{a}(\7h_{a})$. Then $\7h_{a}$ is a
complex Lie subalgebra of $\hol(E,a)$ and $E_{a}\subset E$ is a
complex linear subspace. Denote by $M\subset T$ the subset of all
points at which the function $a\mapsto\dim E_{a}$ takes a global
maximum and fix a connected component $S$ of $M$. Then $S$ is open in
$T$ and $k:=\dim E_{a}$ does not depend on $a\in S$. Let $\2G$ be the
Grassmannian of all $k$-planes in $E$ and consider the map
$\phi:S\to\2G$, $a\mapsto E_{a}$. For every $a\in S$ the map $\phi$ is
constant on $(a+V)\subset S$. Since $\phi$ is CR we conclude that
$T':=E':=E_{a}$ does not depend on $a\in S$. Now fix an arbitrary
vector $\alpha\in E'$ and consider the constant vector field
$\xi=\alpha\dd z$ on $E$. Since $\xi$ is tangent to $S$ it is also
tangent to $T$, that is, the germ $\xi_{a}\in\hol(E,a)$ is contained
in $\7h_{a}$ for all $a\in T$. As a consequence we get $E'\subset
E_{a}$ and thus $E'=E_{a}$ for all $a\in T$. There exists a linear
subspace $V''\subset V$ with $E=E'\oplus E''$ for $E'':=V''\oplus
iV''$. The image $T''$ of $T$ with respect to the canonical projection
$E\to E''$ is a tube submanifold of $E''$ satisfying (iii). The base
$F''$ of $T''$ spans the vector space $V''$, that is, $TÜ''$ is of
finite type by Lemma \ruf{UO}. For the proof of the first part in (ii)
we may assume without loss of generality that $E'=0$ holds, that is,
$E=E''$. But then by the above arguments we have $\7h_{a}=0$ for all
$a\in T$, that is, $T=T''$ is holomorphically nondegenerate.\qed

It is known that for every holomorphically nondegenerate minimal
CR-manifold germ $(M,a)$ the Lie algebra $\hol(M,a)$ has finite
dimension, compare Theorem 12.5.3 in \Lit{BERO}. Calling a CR-manifold
germ $(M, a)$ of {\sl tube type} if it is CR-isomorphic to a germ
$(T,c)$ with $T$ a tube manifold we therefore get the

\Corollary{BZ} Let $(M,a)$ be a CR-manifold germ of tube type. Then
there exist unique integers $k,l\ge0$ and a holomorphically
nondegenerate CR-submanifold $M'\subset M$ of finite type with $a\in
M'$ such that $(M,a)$ is CR-isomorphic to the direct product
$(\CC^{k},0)\times(\RR^{l},0)\times(M',a)$. Furthermore: \0 $(M,a)$ is
holomorphically nondegenerate if and only if $k=0$. \1 $(M,a)$ is of
finite type if and only if $\,l=0$. \1 $\hol(M,a)$ has finite
dimension if and only if $k=l=0$. \Formend\smallskip

As shown in \Lit{BARZ}, to every real-analytic CR-submanifold
$M\subset\CC^{n}$ there exists a proper real-analytic subset $A\subset
M$ such that the germ $(M,a)$ is CR-isomorphic to
$(\CC^{k},0)\times(M',a)$ for some $k\ge0$ and some holomorphically
nondegenerate CR-submanifold $M'\subset M$ containing $a$, provided
$a\in M\backslash A$. Corollary \ruf{BZ} implies that $A$ can be
chosen to be empty if $M$ is of tube type. 

\bigskip \noindent{\bf An analyticity criterion.} In the following
{\sl $k$-differentiable} always means $\5C^{k}$ for $1\le
k\le\infty$. For every abstract $k$-differentiable CR-manifold $N$
then the tangent bundle $\T N$ is of class $\5C^{k-1}$ and we denote
by $\7X^{k-1}(N)$ the $\RR$-linear space of $(k{-}1)$-differentiable
infinitesimal CR-transformations on $N$. Unless $k=k{-}1=\infty$, the
space $\7X^{k-1}(N)$ is not a Lie algebra in general. But again, for
every $k$-differentiable CR-diffeomorphism $\phi:N\to M$ we have a
canonical linear isomorphism
$\phi_{*}:\7X^{k-1}(N)\to\7X^{k-1}(M)$. Clearly, every real-analytic
CR-manifold $M$ can be considered as a $k$-differentiable CR-manifold
in a canonical way and $\hol(M)\subset\7X^{k-1}(M)$ in this sense.

\Proposition{MP} Let $M$ be a real-analytic holomorphically
nondegenerate CR-manifold and let $V+iF$ be a $k$-differentiable tube
submanifold of the complex vector space ${E:=V\oplus\,iV}$. Suppose
that $N$ is an open subset of $V+iF$ and that there exists a
$k$-differentiable CR-diffeomorphism $\phi:N\to M$ with
$\phi_{*}\!\7v\subset\hol(M)$ for $\7v:=\{v\dd z:v\in
V\}\;\subset\;\7X^{k-1}(N)$. Then $N\subset E$ is a (locally-closed)
real-analytic subset of $E$ and $\phi$ is a bianalytic
CR-diffeomorphism.

\Proof Fix an arbitrary point $a\in M$. Since the claim is of local
nature we may assume that $M$ is generically embedded in $E$. The
local flows of vector fields in $\7v$ commute. Therefore the image
$\7w:=\phi_{*}\!\7v$ is an abelian subalgebra of
$\hol(M)\subset\hol(M,a)$ and $\epsilon_{a}(\7w^{\CC})=E$. By
Proposition \ruf{UR} we may assume without loss of generality that
$M=V+iH$ is a real-analytic tube submanifold of $E$ and that
$\7w=\7v\subset\hol(M,a)$. Applying a suitable affine transformation
to $M$ we may assume in addition that $a\in N$, $\phi(a)=a$ and
$\phi:\7v\to\7v$ is the identity. For suitable open subsets
$U,W\subset V$ we may assume furthermore that $F\subset W$, $N=U+iF$
and that there exist $k$-differentiable functions $f,g:U\times W\to V$
satisfying $\phi(z)=f(x,y)+ig(x,y)$ for all $x,y\in U$ with $z=x+iy\in
N$. The condition $\phi_{*}=\id_{\7v}$ implies $\dD fx\equiv\id_{V}$
and $\dD gx\equiv0$ on $U\times F$. The CR-property then gives $\dD
fy|_{c}(v)=0$ and $\dD gy|_{c}(v)=v$ for all $c=(e,f)\in U\times F$
and $v\in \T_{f}F$. Because of $\phi(a)=a$ this implies $\phi(z)=z$ for
all $z\in N$ near $a$, that is, the manifold germs $(N,a)$ and $(M,a)$
coincide.\qed

Proposition \ruf{MP} implies that in case $\7X^{k-1}\!(M,a)=\hol(M,a)$
for every $a\in M$, every $k$-dif\-fer\-en\-tiable tube realization
$N\subset E$ of $M$ is real-analytic. This happens, for instance with
$k=1$, if $M$ is of hypersurface type with nowhere vanishing Levi
form. Indeed, by Theorem 3.1 in \Lit{BAJT} every 1-differentiable
CR-diffeomorphism between open subsets of $M$ is real-analytic.

\KAP{realizations}{Tube realizations}

In the following $M$ is a CR-manifold generically embedded in the
complex manifold $Z$ and $a\in M$ is a given point. Then the tube
realizations $\phi:(M,a)\to(T,c)$ and $\phi:(M,a)\to(T',c')$ with
tubes $T\subset E$, $T'\subset E'$ as in Section \ruf{manifolds} are
called {\sl affinely equivalent} if the tube germs $(T,c)$, $(T',c')$
are equivalent under an affine isomorphism $\lambda:E\to E'$, or
equivalently, if $\phi'\circ g=\lambda\circ\phi$ for some
$g\in\Aut(M,a)$ and some affine isomorphism $\lambda$. Also, we call
the subsets $\7v,\7v'\subset\hol(M,a)$ conjugate with respect to
$\Aut(M,a)$ if $\7v'=g_{*}(\7v)$ for some $g\in\Aut(M,a)$.

\Proposition{UR} The affine equivalence classes of tube realizations of
the germ $(M,a)$ are in 1-1-cor\-re\-spon\-dence to the 
$\Aut(M,a)$-conjugacy classes of abelian Lie subalgebras
$\7v\subset\7g:=\hol(M,a)$ satisfying \0 $\7v$ is totally real in
$\hol(Z,a)$, \1 $\7e:=\7v\oplus\,i\!\7v\subset\hol(Z,a)$ spans the full
tangent space $\T_{a}Z$.

\Proof Suppose that for the tube submanifold $T=V+iF\subset E$ the
CR-isomorphism $\phi:(M,a)\to(T,c)$ is given. Then
$\7v:=\{\phi^{*}(v\dd z):v\in V\}\subset\7g$ satisfies the
conditions (i) and (ii). For every affine isomorphism $\lambda:E\to
E'$ and every tube realization $\phi':(M,a)\to(T',c'):=\lambda(T,c)$
the transformation $g:=\phi'^{-1}\circ\lambda\circ\phi\in\Aut(M,a)$
satisfies $g_{*}(\7v)=\{\phi'^{*}(v'\dd z):v'\in V'\}$ for
$V':=\lambda(V)$.\nline Conversely, suppose that an abelian Lie
subalgebra $\7v\subset\7g$ with (i), (ii) is given.  Then $\7e$
is an abelian complex Lie algebra and by Lemma \ruf{LU} the evaluation
map $\epsilon_{a}:\7e\to\T_{a}Z$ is a complex linear
isomorphism. Denote by $E$ the complex vector space underlying $\7e$
and by $V\subset E$ the real vector space underlying $\7v$. By the
implicit function theorem there exist open neighbourhoods $U$ of $0\in
E$ and $W$ of $a\in Z$ such that $\psi(\xi):=\exp(\xi)(a)\in W$ is
defined for every vector field $\xi\in U$ and $\psi:U\to W$ is a
biholomorphic mapping with $\psi(0)=a$. For $\phi:=\psi^{-1}$ then
$\phi(W\cap M)$ is an open piece of a tube $T=V+iF\subset E$, that is,
$\phi:(M,a)\to(T,0)$ gives a tube realization with
$\7v=\{\phi^{*}(v\dd z):v\in V\}$. Now fix a $g\in\Aut(M,a)$. Then
also $\7v':=g_{*}(\7v)$ with $\7e':=\7v'\oplus\,i\!\7v'$ satisfies
(i), (ii) and thus gives a tube realization $\phi':(M,a)\to(T',0)$
according to the procedure above. Since $\7e$, $\7e'$ are abelian,
there is a complex linear isomorphism $\lambda:E\to E'$ with
$\lambda(\psi(\xi))=\psi'(g_{*}(\xi))$ for all $\xi$ in a
neighbourhood of the origin in $E$. But this means that
$\lambda:(T,0)\to(T',0)$ is an affine equivalence.\qed

Notice that $\7e$ is maximal abelian in $\hol(Z,a)$ by Lemma
\ruf{LU}. In case $M$ is holomorphically non-degenerate the condition
(i) above is automatically satisfied and $\7v$ is maximal abelian in
$\hol(M,a)$.

\Remark{} A different characterization of abelian Lie subalgebras
$\7v$ giving rise to tube realizations of $(M,a)$ occurs already in
\Lit{BART}. Instead of (i), (ii) there $\7v$ has to act without
isotropy and transversally to the holomorphic tangent bundle.

\bigskip Tubes $T=V+iF$ have a special property: $\tau(x+iy):=-x+iy$
for all $x\in V$, $y\in F$ defines an anti-CR map $\tau:T\to T$ with
$\tau^{2}=\id$ and $\tau(a)=a$ for all $a\in iF\subset T$. This
motivates the following considerations.

{\bf\smallskip\noindent Involutions.}  In this subsection $M$ stands
 for an arbitrary CR-manifold. A real-analytic mapping $\tau:M\to M$
 is called an {\sl involution} of $M$ if it is anti-CR and satisfies
 $\tau^{2}=\id$. If in addition $\tau(a)=a$ for a given $a\in M$ we
 call $\tau$ an {\sl involution of $M$ at} $a$ or of the CR-manifold
 germ $(M,a)$. Two involutions $\tau$, $\tau'$ of $(M,a)$ are called
 {\sl equivalent} if $\tau'=g\tau g^{-1}$ for some
 $g\in\Aut(M,a)$. Every involution $\tau$ of $(M,a)$ splits various
 linear spaces, associated with the germ $(M,a)$, into their
 $\pm1$-eigenspaces. To indicate the dependence on $\tau$ we mark the
 $+1$-eigenspaces by an upper index $\tau$ and the $-1$-eigenspaces by
 an upper index $-\tau$, e.g. $\T_{a}M=\T^{\tau}_{a}M\oplus
 \T^{-\tau}_{a}M$, $\H_{a}M=\H^{\tau}_{a}M\oplus \H^{-\tau}_{a}M$ and
 $\7g=\7g^{\tau}\oplus\7g^{-\tau}$ for $\7g:=\hol(M,a)$. Clearly
 $(\T_{\!a}M)^\tau=\T_{\!a}(M^{\tau})$. 
Crucial for the explicit determination of all tube realizations for
 $(M,a)$ is the following reformulation of Proposition \ruf{UR},
 compare e.g. \Lit{FLKA}.

\Proposition{KR} Proposition \ruf{UR} remains valid if $\,${\rm (i)}
is replaced by \item{\rm (i')} There exists an involution $\tau$ of
$(M,a)$ with $\7v\subset\7g^{-\tau}$.\par\vskip2pt \noindent The
involution $\tau$ in {\rm (i')} is uniquely determined by $\7v$ and
satisfies \1 $\dim_{a}\!M^{\tau}=\hbox{\rm CR-}\!\dim M$, or
equivalently, $\,\dim\T^{-\tau}_{a}M=\dim_{\CC}\!Z$.\par\vskip2pt
\noindent In particular, for every $g\in\Aut(M,a)$ the involution
$\tau'$ corresponding to $\7v':=g_{*}(\7v)\subset\7g$ is given by
$\tau'=g\circ\tau\circ g^{-1}$.

\Proof (i') $\Longrightarrow$ (i) is obvious. Therefore let us assume
conversely that the abelian subalgebra $\7v\subset\7g$ satisfies (i)
and (ii). Without loss of generality we assume by Proposition \ruf{UR}
that $M=V+iF$ is a tube submanifold of $E=V\oplus iV$, that $a\in iF$
and that $\7v=\{v\dd z:v\in V\}$. Then the involution
$\tau(x+iy)=-x+iy$ of $(M,a)$ satisfies (i') and (iii). Now suppose
that $\tau'$ is a further involution of $(M,a)$ with the same
properties. Then $g:=\tau\circ\tau'\in\Aut(M,a)$ satisfies
$g_{*}(\alpha\dd z)=\alpha\dd z$ for all $\alpha\in V$ and hence also
for all $\alpha\in E$. But then $g=\id$ and $\tau'=\tau$ .\qed

\Remark{} The explicit determination of all tube realizations for
$(M,a)$ up to affine equivalence requires by Proposition \ruf{UR}
that, up to conjugation by the stability group $\Aut(M,a)$, all
abelian Lie subalgebras $\7v\subset\hol(M,a)$ have to be found that
satisfy the conditions (i), (ii). Proposition \ruf{KR} restricts the
search (and with it the amount of computation) to the following:
Determine first, up to conjugation, all involutions of $(M,a)$ that
satisfy (iii) and then, for every such involution $\tau$, search
for suitable $\7v$'s in $\7g^{-\tau}$.  As an application of that
method we classify algebraically in the forthcoming paper \Lit{FLKA}
all local tube realizations of Levi nondegenerate hyperquadrics
$Q\subset\CC^{n}$.  These are locally CR-equivalent to the
hypersurfaces $S^{1}_{pq}$ occurring in the next section and have the
special property that every germ $(Q,a)$ has, up to equivalence, a
unique involution satisfying (iii).

\KAP{Classification}{Classification of involutions for a special class
of CR-manifolds}

Fix in the following arbitrary integers $p,q\ge m\ge1$ and denote by
$\2G$ the Grassmannian of all linear $m$-spaces in $\CC^{n}$,
$n:=p+q$.  Then $\2G$ is a compact complex manifold of dimension
$m(n-m)$ on which $\SL(n,\CC)$ acts transitively by holomorphic
transformations.  The group $\Aut(\2G)$ coincides with
$\PSL(n,\CC)=\Quot{\SL(n,\CC)}{\hbox{center}}$, unless $p=q>1$ (in
which case there is a second connected component of $\Aut(\2G)$). To
avoid totally real examples we exclude the case $p=q=m$ for the rest
of the section.

Consider on $\CC^{n}$ the real-valued function $\6h$ defined by
$$ \6h(z)=(u|u)-(v|v)\Steil{for all}z=(u,v)\in\CC^{p}\oplus\CC^{q}$$
with $(\,|\,)$ being the standard inner product and identify
$\SU(p,q)\subset\SL(n,\CC)$ with the subgroup of all transformations
leaving $\6h$ invariant. Then the connected real submanifold
$$S:=S_{m}^{p,q}:=\big\{L\in\2G:\6h(L)=0\big\}\Leqno{TU}$$ is the
unique closed (and hence compact) $\SU(p,q)$-orbit in $\2G$. As
CR-submanifold $S$ is generically embedded in $\2G$ with CR-dimension
$m(n-2m)$ and CR-codimension $m^{2}$. Furthermore, a dense open subset
of $S$ can be realized as a real quadric in $\CC^{n-m}$,
$\7g:=\su(p,q)=\hol(S)\cong\hol(S,a)$ holds for every $a\in S$,
compare \Lit{KAPW} for details. As a consequence of Theorem 1.3 in
\Lit{ISKA} every CR-isomorphism between domains $D_{1},D_{2}$ of $S$
extends to a biholomorphic automorphism of $\2G$ leaving $S$
invariant. Since $S$ has a global (anti-CR) involution (see the
following classification) also every anti-CR-isomorphism between
domains $D_{1},D_{2}$ of $S$ extends to a global antiholomorphic
automorphism of $\2G$ leaving $S$ invariant. For the classification of
all involutions of the germ $(S,a)$ it is therefore enough to
determine all {\sl global} involutions of $S$.

\medskip\noindent {\bf Classification of all involutions on S.}  Let a
global involution $\tau$ of $S$ (not necessarily having a fixed point)
be given. Then $\tau$ extends to an antiholomorphic automorphism of
$\2G$ that we also denote by $\tau$. Also, the involution induced by $\tau$
on $\7l:=\7{sl}(n,\CC)\cong\aut(\2G)$ will be denoted by the same
symbol. The fixed point submanifold $\2G^{\tau}$ of $\2G$ is either empty or a
real form of $\2G$. One can show that there are integers
$\epsilon,\delta$ with $\epsilon^{2}=\delta^{2}=1$ together with an
antilinear endomorphism $\tild\tau$ of $\CC^{n}$ such that
$\tild\tau^{2}=\epsilon\id$, $\,\6h\circ\tild\tau=\delta\6h$ and
$\tau(L)=\{\tild\tau(z):z\in L\}$ for all $L\in\2G$. Depending on the
value of $\epsilon$ we have the following two cases.
{\advance\parindent17pt \item{\hfill$\epsilon=\phantom{-}1:$} Then
$\7l^{\tau}\cong\7{sl}(n,\RR)$ and $\2G^{\tau}$ can be identified with
the real Grassmannian of all real linear m-spaces in $\RR^{n}$.
\item{\hfill$\epsilon=-1:$} This case can only occur if $n$ is even
and then $\7l^{\tau}\cong\7{sl}(n/2,\HH)$, where $\HH$ is the field of
quaternions. Furthermore, $\2G^{\tau}$ is empty if and only if $m$ is
odd.\par }

\medskip\noindent
The precise classification requires some work. Here we state only the
final result: It turns out that for every given $p,q$ the possible
pairs $(\epsilon,\delta)$ stand in a one-to-one relation with the
$\SU(p,q)$-conjugation classes of involutions on $S=S_{m}^{p,q}$. More
explicitly,  every such involution is conjugate to exactly one of
the following types {\bf I} -- {\bf IV}, where we write every (row)
$z\in\CC^{n}$ in the form $z=(u,v)$ with $u\in\CC^{p}$ and
$v\in\CC^{q}$. Also, for every integer $d\ge1$ we put
$J_{d}:=\pmatrix{0&\One\cr-\One&0\cr}\in\GL(2d,\2Z)\,.$

\item{\bf I:} $(\epsilon,\delta) =(1,1)$ and $\tild\tau(z)=\overline
z$. The fixed point set $S^{\tau}$ has dimension $m(n-2m)$ and is an
orbit of the subgroup $\SO(p,q)\subset\SU(p,q)$. Also $\7l^\tau=
\7{sl}(n,\RR)$ and $ \7g^\tau=\so(p,q)$ for the $\tau$-fixed point
subsets.

\item{\bf II:} $(\epsilon,\delta) =(1,-1)$, $\;p=q$ and
$\tild\tau(z)=(\overline v,\overline u)$. Here $S^{\tau}=\2G^{\tau}$
has dimension $m(n-m)$. Also $\7l^\tau\cong \7{sl}(2p,\RR)$ and
$\7g^\tau\cong \7{sp}(p,\RR)$.

\item{\bf III:} $(\epsilon,\delta) =(-1,1)$, $\;p=2p'$, $q=2q'$ are
even and $\tild\tau(z)=(\overline uJ_{p'},\overline vJ_{q'})$. Here
$\7l^\tau\cong \7{sl}(p'+q',\HH)$ and $\7g^\tau\cong \7{sp}(p',q')$.

\item{\bf IV:} $(\epsilon,\delta) =(-1,-1)$, $\;p=q$ and
$\tild\tau(z)=\overline zJ_{p}=(-\overline v,\overline u)$. Here
$\7l^\tau\cong \7{sl}(p,\HH)$ and $\7g^\tau\cong \so(p,\HH)$.
Furthermore, $S^{\tau}$ has dimension $m(n-m-1)$ if $m$ is even.

\bigskip\noindent From the above classification we see that in case
$m$ odd for every $a\in S$ there exists exactly one conjugation class
of involutions of the germ $(S,a)$ which satisfies condition (iii) in
Proposition \ruf{KR}
(namely the one given by type {\bf I} above). Also, for the types {\bf
II}, {\bf IV}, every $m$ and every $a\in S^{\tau}$ the germ $(S,a)$
never satisfies condition  (iii) in \ruf{KR}.

\KAP{coarser}{A coarser equivalence relation}

In general, for a given tube submanifold $T=V+iF$ of $E=V+iV$, there
is an infinite subset $A\subset T$ such that for every $a\ne b$ in $A$
the germs $(T,a)$, $(T,b)$ are affinely inequivalent - even if $T$ is
locally homogeneous and hence all $(T,a)$, $(T,b)$ are
CR-equivalent. As an example consider in $\CC^{2}$ the closed tube
hypersurface $T=\RR^{2}+iF$ with
$$F:=\{x\in\RR^{2}:\cos x_{1}=e^{x_{2}},\,|x_{1}|<\pi/2\}$$ (the
boundary of the middle gray domain in Figure 1, Section
\ruf{sphere}). Consider the function $f(z):=\Im(z_{2})$ on $T$. Then
for every $a,b\in T$ the germs $(T,a)$, $(T,b)$ are CR-isomorphic, (in
fact, $T$ is locally CR-isomorphic to the euclidian sphere
$S^{3}\subset\CC^{2}$) but they are affinely equivalent if and only if
$f(a)=f(b)$. Therefore $T$ gives rise to a continuum of mutually
affinely inequivalent tube realizations of the CR-germ $(T,0)$.  This
phenomenon motivates the introduction of a coarser equivalence
relation that puts all germs $(T,a)$, $a\in T$, into a single
equivalence class. The construction is motivated by the concept of a
sheaf:\nline For fixed $E=V^{\CC}$ let $\5T=\5T(V)$ be the set of all
(real-analytic) germs $(T,a)$ with $T=V+iF$ an arbitrary tube
submanifold of $E$ and $a\in T$. Furthermore define $\pi:\5T\to E$ by
$(T,a)\mapsto a$. Then $\5T$ becomes in the standard way a Hausdorff
topological space over $E$ -- the topology on $\5T$ is the coarsest
one such that for every tube submanifold $T\subset E$ the subset
$[T]:=\{(T,a):a\in T\}$ is open in $\5T$. The space $\5T$ has in a
unique way the structure of a (disconnected) CR-manifold by requiring
that $\pi:[T]\to T$ is a CR-isomorphism for every tube submanifold
$T\subset E$. Every real affine transformation
$g\in\Aff(V)\subset\Aff(E)$ (the respective affine transformation
groups) gives rise to a CR-automorphism of $\5T$
by $g(T,a):=(gT,ga)$, that we also denote by $g$. However, it should
be noticed that the corresponding action of the Lie group $\Aff(V)$ on
$\5T$ is discontinuous. Nevertheless, every connected component of
$\5T$ is invariant under the (continuous) action of the translation
subgroup $V\subset\Aff(V)$ and therefore may be considered as a {\sl
generalized tube manifold over $E$}. For every (connected) tube
submanifold $T\subset E$ denote by $\tilde T$ the connected component
of $\5T$ containing $[T]$ and call the pair $(\tilde T,\pi)$ the {\sl
abstract globalization} of $T$ and also of every tube germ $(T,a)$,
$a\in T$.  Since the translation group $V\subset\Aff(E)$ acts on
$\tilde T$ by CR-transformations we may consider $\tilde T$ as tube
manifold {\sl over} $E$ via $\pi$.

\Definition{PF} The tube manifold germs $(T,a)$, $(T',a')$ in
$E=V^{\CC}$ are called {\sl globally affinely equivalent} if
$\;\tilde{T'}=g(\tilde T)$ for the corresponding abstract globalizations
and a suitable $g\in\Aff(V)$.

In case $\pi(\tilde T)$ is a (locally closed) submanifold of $E$, we
call $\pi(\tilde T)$ the {\sl globalization} of $(T,a)$ and denote it
by $\hat T$. Clearly, then $\pi:\tilde T\to\hat T$ is a
CR-isomorphism. Furthermore, $\hat T$ is a tube submanifold of $E$
containing $T$ as an open submanifold and also is maximal with respect
to this property. As an example, every closed tube submanifold
$T\subset E$ is the globalization of each of its germs $(T,a)$, $a\in
T$.

\medskip In the following we assume for the CR-manifold germ $(M,a)$
that the Lie algebra $\7g:=\hol(M,a)$ has finite dimension. Then, in
particular, $(M,a)$ is holomorphically nondegenerate and we denote as
usual with $\Int(\7g)\subset\Aut(\7g)$ the inner automorphism group of
$\7g$, that is, the subgroup generated by all $\exp(\ad\xi)$,
$\xi\in\7g$. Finally, for every $a\in M$ let
\par \hfil $\rho_{a}:\hol(M)\hookrightarrow\hol(M,a)$ \hfil
\smallskip\noindent be the restriction mapping. Then we have

\Lemma{LE} Suppose that $\rho_{a}:\hol(M)\to\7g=\hol(M,a)$ is
bijective. Then $g\mapsto\rho_{a}g_{*}\rho_{a}^{-1}$ defines a group
homomorphism $\Aut(M)\to\Aut(\7g)$ that sends $H$ to $\Int(\7g)$,
where $H\subset\Aut(M)$ denotes the subgroup generated by
$\exp(\aut(M))$. For every $g\in\Aut(M)$ and $b:=g(a)$ also
$\rho_{b}\!:\hol(M)\to\hol(M,b)$ is bijective. Furthermore, for every
abelian subalgebra $\7w\subset\hol(M)$ such that
$\rho_{a}(\7w)\subset\7g$ gives a local tube realization, also
$\rho_{b}(g_{*}\!\7w)\subset\hol(M,b)$ gives a local tube realization
and both are affinely equivalent.

\Proof From $\rho_{a}\!\Ad(\exp\xi)=\rho_{a}\!\exp(\ad\xi)=
\exp(\ad\rho_{a}(\xi))\rho_{a}$ for all $\xi\in\aut(M)$ we see that $H$
maps into $\Int(\7g)$. The other statements are obvious.\qed

\medskip The following global statement will be one of the key
ingredients in the proof of the following Theorem \ruf{HT}. Both of
these results allow to reduce the classification problem for tube
realizations of $(M,a)$ in many cases to a purely algebraic one.

\Proposition{ZB} Let $Z$ be a complex manifold and $M\subset Z$ a
generically embedded minimal CR-submani\-fold. Assume that, for a given
point $a\in M$, $\7g:=\hol(M,a)$ has finite dimension and every germ
in $\7g$ extends to a vector field in $\aut(M)$. Let
$\7v,\7v'\subset\7g$ be abelian subalgebras giving rise to local tube
realizations of $(M,a)$ according to Proposition \ruf{UR} and assume
that every germ in $\7e:=\7v^{\CC}\subset\hol(Z,a)$ extends
to a vector field in $\aut(Z)$. Then the local tube realizations of
$(M,a)$ given by $\7v,\7v'$ are globally affinely equivalent if
$\7v=\lambda(\7v')$ for some $\lambda\in\Int(\7g)$.

\Proof For simpler notation we identify the Lie algebras $\hol(M)$ and
$\7g $ via the isomorphism $\rho_{a}:\hol(M)\to\7g$. Since
$\aut(M)=\hol(M)$, for every $\lambda\in\Int(\7g)$ with
$\7v=\lambda(\7v')$ there exists a transformation $g\in G$ with
$\lambda=\Ad(g)=g_{*}$, where $G$ is the group $H$ from Lemma
\ruf{LE}. For $b:=g(a)$ the abelian subalgebras $\7v'\subset\7g$ and
$\rho_{b}(\7v)\subset\hol(M,b)$ give affinely equivalent local tube
realizations. Therefore we have to show that the abelian subalgebras
$\7v\subset\7g$ and $\rho_{b}(\7v)\subset\hol(M,b)$ give globally
affinely equivalent tube realizations of the germs $(M,a)$ and
$(M,b)$.\quad To begin with let $E$ and $V$ be the vector spaces
underlying $\7e$ and $\7v$, compare the proof of Proposition
\ruf{KR}. Then the locally biholomorphic map $\psi:E\to Z$,
$\xi\mapsto\exp(\xi)(a)$, is the universal covering of an open subset
$O\subset Z$ with $Z\backslash O$ analytic in $Z$. Denote by $T$ the
connected component of $\psi^{-1}(M)$ that contains the origin of
$E$. By Lemma \ruf{BR} the intersection $M\cap\,O$ is connected, that
is, there is a point $c\in T$ with $\psi(c)=b$.  Now $T$ is a tube
submanifold of $E$ and the tube germ $(T,0)$ is affinely equivalent to
the tube realization of $(M,a)$ given by $\7v\subset\hol(M,a)$. Also
the tube germ $(T,c)$ is affinely equivalent to the tube realization
of $(M,b)$ given by $\rho_{b}(\7v)\subset\hol(M,b)$. This proves the
claim.\qed

\Corollary{CO} In case $M$ in Proposition \ruf{ZB} is closed in $Z$,
the tube realization of $(M,a)$ given by $\7v$ is affinely equivalent
to the germ $(T,0)$ with $T\subset E$ a suitable $\,${\rm closed}$\,$
tube submanifold containing the origin. In other words, the germ
$(T,0)$ has a closed globalization $\hat T$ in $E$.

\Proof With the notation of the proof for Proposition \ruf{ZB} the
intersection $M\cap O$ is closed in $O$. Hence also $T\subset E$ is
closed.\qed

Since every $M=S^{m}_{pq}$, compare \Ruf{TU}, is closed in $Z=\2G$ and
the assumptions of Proposition \ruf{ZB} are satisfied for $M\subset
Z$, we have: {\sl Every tube submanifold of $\CC^{r}$ locally
CR-equivalent to $S^{m}_{pq}$ extends to a closed tube submanifold of
$\CC^{r}$ with the same property}. For the special case $m=1$ this
statement is already contained in \Lit{ISAV}.  \nline In case the
manifold $M$ is not closed in $Z$ the globalization of a local tube
realization for $M$ may be no longer closed in $E$. For a typical
example compare the lines following \Ruf{ST}.

\KAP{The}{The subgroup $\,\Glob(M,a)\subset\Aut(\hol(M,a))$}

\noindent In certain cases also the converse of Proposition \ruf{ZB}
is true. Let us denote for $\7g=\hol(M,a)$ by
$$\Glob(M,a)\subset\Aut(\7g)\hbox{ the subgroup generated by }$$
$$\Int(\7g)\Steil{together
with}\Ad\big(\Aut(M,a)\big)=\{g_{*}:g\in\Aut(M,a)\}\,.$$ Clearly,
$\Int(\7g)$ is a connected subgroup of $\Glob(M,a)$ and coincides with
the connected identity component of $\Aut(\7g)$ if $\7g$ is
semi-simple. For the complex manifold $Z$ and the CR-submanifold
$M\subset Z$ we will need the following

\noindent {\bf Condition P:} {\sl Every CR-isomorphism of germs
$(M,a)\to(M,b)$ with $a,b\in M$ extends to an automorphism
$g\in\Aut(Z)$ with $g(M)=M$.}

\noindent {\bf Condition Q:} There exists an antiholomorphic
automorphism $\tau$ of $Z$ with $\tau(M)=M$.

Notice that if Conditions {\rm P} and {\rm Q} are satisfied for
$M\subset Z$ simultaneously then also every anti-CR-isomorphism of
germs $\theta:(M,a)\to(M,b)$, $a,b\in M$, extends to an
antiholomorphic automorphism $\theta$ of $Z$ leaving $M$
invariant. Indeed, for $c:=\tau(b)$ the CR-isomorphism
$\tau\circ\theta:(M,a)\to(M,c)$ extends to a $g\in\Aut(Z)$ with
$g(M)=M$. But then $\tau^{-1}\circ g$ is the antiholomorphic extension
of $\theta$ to $Z$.

\Theorem{HT} Let $Z$ be a compact complex manifold and $M\subset Z$ a
homogeneous generically embedded closed CR-submanifold satisfying
condition {\rm P}. Then, given $a\in M$, any two local tube
realizations of the germ $(M,a)$ given by the abelian Lie subalgebras
$\7v,\7v'\subset\7g$ are globally affinely equivalent if and only if
$\7v=\lambda(\7v')$ for some $\lambda\in\Glob(M,a)$.
Furthermore, the Lie algebra
$\7g:=\hol(M,a)$ has finite dimension.

\Proof $\Aut(Z)$ is a complex Lie group in the compact-open topology
with Lie algebra $\aut(Z)=\hol(Z)$ since $Z$ is compact. Every
$\xi\in\7g$ defines a local flow in a small open neighbourhood of
$a\in M$ and thus a one parameter subgroup of $\Aut(Z)$ by condition
{\rm P}. Therefore every such $\xi$ extends to a vector field in
$\hol(Z)$ tangent to $M$. Identifying $\7g$ and $\hol(M)$ as before
via the isomorphism $\rho_{a}$ we have $\7g=\hol(M)\subset\hol(Z)$. In
particular, $\7g$ has finite dimension. Let $G\subset\Aut(M)$ be the
subgroup generated by $\exp(\aut(M))$. Then $G$ acts transitively on
$M$ since by assumption $M$ is homogeneous. Therefore every
$g\in\Aut(M)$ is of the form $g=g_{1}g_{2}$ with $g_{1}\in G$ and
$g_{2}(a)=a$. This implies
$$\Ad(\Aut(M))\;\subset\;\Int(\7g)\Ad(\Aut(M,a))\;=
\;\Glob(M,a)\,.\leqno{(*)}$$ {\bf `if'} In case
$\lambda\in\Ad(\Aut(M,a))$ the abelian Lie algebras $\7v,\7v'$ already
give affine equivalent local tube realizations of $(M,a)$ by
Proposition \ruf{UR}. It is therefore enough to discuss the case
$\lambda\in\Int(\7g)$. But this follows immediately with Proposition
\ruf{ZB}.  \nline {\bf `only if'} By Corollary \ruf{CO} there are
closed tube submanifolds $T,T'$ of $E=V^{\CC}$ containing the origin
such that $(T,0)$ and $(T',0)$ are the local tube realizations of
$(M,a)$ determined by $\7v$ and $\7v'$. Also there are locally
biholomorphic mappings $\psi,\psi':E\to Z$ with $\psi(0)=\psi'(0)=a$
and such that $\psi(T)$ as well as $\psi'(T')$ are open in $M$,
compare the proof of Proposition \ruf{ZB}. Now assume that $(T,0)$ and
$(T',0)$ are globally affinely equivalent. Then there exists a complex
affine automorphism $h$ of $E$ with $T=h(T')$ (but not necessarily
with $h(0)=0$).  By condition {\rm P} there is a unique $g\in\Aut(Z)$
with \nline\line{$(\dagger)$\hss$g\circ\psi'=\psi\circ
h$\hss}\vskip0pt\noindent and $g(M)=M$.  Put $b:=g(a)$ and
$c:=h(0)$. Then $\psi(c)=b$ and $\lambda:=\Ad(g)\in\Glob(M,a)$ by
$(*)$. For $\7V:=\{v\dd z:v\in V\}\subset\hol(E)$ we have\nline
\line{\hss$h_{*}(\rho_{0}(\7V))=\rho_{c}(\7V)$,\quad
$\psi_{*}(\rho_{c}(\7V))=\rho_{b}(\7v)$\steil{and}
$\psi'_{*}(\rho_{0}(\7V))=\rho_{a}(\7v')$,\hss} where $\rho_b$ is the
restriction map, introduced just before \ruf{LE}.  This implies
$\rho_{b}(\7v)=\rho_{b}(\lambda(\7v'))$ by $(\dagger)$ and hence
$\7v=\lambda(\7v')$ as desired. \qed

\medskip{\noindent\bf An example for Theorem \HT.} As an example for
a pair $M\subset Z$ satisfying all the assumptions in \ruf{HT} we may
take the complex projective space $Z=\PP_{r}$ together with the
compact homogeneous hypersurface $S=S_{1}^{p,q}$ from \Ruf{TU} as $M$,
where the integers $p,q,r\ge1$ satisfy $p+q=r+1\ge3$. Condition {\rm P}
for example, is satisfied by Theorem 6 in \Lit{TANA}. Then
$L:=\Aut(Z)=\PSL(r+1,\CC)$ and $G:=\{g\in L:g(S)=S\}$ can be
canonically identified with $\Aut(M)$. The real Lie group $G$ has
$(1+\delta_{p,q})$ connected components, the connected identity
component $G^{0}=\PSU(p,q)$ is a real form of $L^{0}$. For the Lie
algebras we have $\7l:=\hol(Z)=\7{sl}(r+1,\CC)$ with real form
$\7g:=\hol(S)=\su(p,q)$. If we fix $a\in S$ and identify the Lie
algebras $\7g$, $\hol(S,a)$ via the restriction operator $\rho_{a}$ we
have $\Glob(S,a)=\Ad(G)\cong G$. In particular, $\Glob(S,a)=\Int(\7g)$
if $p\ne q$.

Now suppose that $\7e\subset\7l$ is a complex abelian subalgebra such
that the subgroup $\exp(\7e)\subset L$ has an open orbit $O$ in
$\PP_{r}$. By Lemma \ruf{LU} then $\7e$ has dimension $r$ and is
maximal abelian in $\7l$. The orbit $O$ consists of all points $c\in
Z=\PP_{r}$ with $\epsilon_{c}(\7e)=\T_{c}Z$, and the complement
$A:=\PP_{r}\!\backslash O$ is the union $A=H_{1}\cup
H_{2}\cup\cdots\cup H_{k}$ of $k\le r+1$ complex projective
hyperplanes $H_{j}$ in $\PP_{r}$. Clearly, the conjugacy class of
$\7e$ in $\7l$ modulo the action of $L$ depends on the $L$-orbit of
$A$ in the space of all analytic subsets in $\PP_{r}$.

Suppose in addition that $\7e=\7v^{\CC}$ for $\7v:=\7e\cap\7g$ and fix
a point $a\in O\cap S$. Put $E:=\CC^{r}$, $V:=i\RR^{r}$ and choose a
complex linear isomorphism $\Xi:E\to\7e$ with $\Xi(V)=\7v$. Then the
locally biholomorphic map $\psi:E\to O$, $\,z\mapsto\exp(\Xi z)a,\,$
realizes $E=\CC^{r}$ as universal cover of the domain $O$. The
intersection $O\cap S$ is a closed CR-submanifold of $O$ and divides
$O\backslash S$ into the two connected components $O^{\pm}:=O\cap
S^{\pm}$. In general the pre-image $\psi^{-1}(S)$ in $E$ decomposes
into several connected components which only differ by a translation
in $E$. Let $T$ be one of these. Then by Corollary \ruf{CO} $T$ is a
closed tube submanifold of $E$ and a covering of $O\cap S$ via $\phi$.

In the next section we will discuss the special case $p=1$.

\KAP{sphere}{The standard sphere}

\noindent In this section we consider for fixed $r\ge2$ the euclidian
hypersphere
$$S:=\{z\in\CC^{r}:(z|z)=\sum z_{k}\overline z_{k}=1\}\,.\Leqno{AR}$$
$S$ is the boundary of the euclidian ball
$B:=\{z\in\CC^{r}:(z|z)<1\}$, a bounded symmetric domain of rank 1. We
always consider $\CC^{r}$ as domain in the complex projective space
$\PP_{r}$ by identifying the points $(z_{1},\dots,z_{r})\in\CC^{r}$
and $[1,z_{1},\cdots,z_{r}]\in\PP_{r}$. In this sense $S$ can also be
written as $$S=\big\{[z_{0},\cdots,z_{r}]\in\PP_{r}:z_{0}\overline
z_{0}=\sum_{k>0}z_{k}\overline z_{k}\big\}\,,$$ which is the case
$p=1,$ $q=r$ at the end of the preceding section. Every $g\in\Aut(S)$
extends to a biholomorphic automorphism of $\PP_{r}$ leaving the ball
$B=S^{+}$ as well as the outer domain $\PP_{r}\!\backslash\overline
B=S^{-}$ invariant and thus gives isomorphisms of the groups
$$G:=\Aut(S)\cong\Aut(S^{\pm})\cong\{g\in\Aut(\PP_{r}):g(S)=S\}
\cong\PSU(1,r)\,,$$ which we identify in the following. In particular,
$S$ is homogeneous and $G$ is a real form of
$L:=\Aut(\PP_{r})=\PSL(r+1,\CC)$. It is well known that
$\Aut_{a}(S)=\Aut(S,a)$ holds for every $a\in S$ and that
$\Aut(S,a)$ acts transitively on the ball $B$.

In the following we describe some abelian Lie subalgebras
$\7v\subset\7g:=\hol(S)$ that lead to local tube realizations of
$S$. Every vector field in $\7l:=\7g^{\CC}=\hol(\PP_{r})$ is
polynomial of degree $\le2$ in the coordinate $z=(z_{1},\dots,z_{r})$
of $\CC^{r}$ and
$$\7g=\big\{(\alpha+zu-(z|\alpha)z)\dd
z:\alpha\in\CC^{r},u\in\7u(r)\big\}\,.$$ With $E:=\CC^{r}$ and
$V:=i\RR^{r}$ we start with an arbitrary but fixed $\alpha\in V$ and
consider the abelian subalgebra
$$\7v:=\RR\big(\alpha-(z|\alpha)z\big)\dd z\,\oplus\;\{zu\dd
z:u\in\7u(r)\steil{diagonal with}\alpha u=0\}\subset\7g\,.$$ Then
$\7e:=\7v^{\CC}\subset\7l$ has an open orbit $O\subset\PP_{r}$ and,
fixing a complex linear isomorphism $\Xi:\CC^{r}\to\7e$ as at the end
of the preceding section, we get the universal covering map
$\phi:\CC^{r}\to O$.\nline {\bf In case $\alpha=0$} we have
$O=(\CC^{*})^{r}$ and $\phi$ can be chosen as
$\phi(z)=(e^{z_{1}},\dots,e^{z_{r}})$. Then
$T:=\phi^{-1}(S)=F+i\RR^{r}$ is the tube with base
$$F=\{x\in\RR^{r}:e^{2x_{1}}+e^{2x_{2}}+\dots+e^{2x_{r}} =1\}\,.$$
With $e^{2x_{1}}-1=2e^{x_{1}}\sinh x_{1}$ it is obvious that $F$ is
affinely equivalent in $\RR^{r}$ to the hypersurface
$$\Pi_{-}:=\big\{x\in\RR^{r}:\sinh x_{1}=\sum_{k>1}e^{x_{k}}\big\}$$
occurring in Theorem 2 of \Lit{DAYA}. Notice that $\PP_{r}\backslash
O=H_{0}\cup H_{1}\cup\cdots\cup H_{r}$ is the union of $r{+}1$
projective hyperplanes in general position with $H_{1},\dots,H_{r}$
intersecting $S$ transversally and $H_{0}$ not meeting $S$.

\noindent {\bf In case $\alpha=(i,0,\dots,0)$} we have
$O=\big\{[z]\in\PP_{r}:(z_{0}^{2}+z_{1}^{2})z_{2}z_{3}\cdots
z_{r}\ne0\big\}$
and $\phi$ can be chosen as 
$$\phi(z)=[\cos z_{1},\sin z_{1},e^{z_{2}},\dots,e^{z_{r}}]\steil{for
all}z\in\CC^{r}\,.$$ $\phi^{-1}(S)$ has a countable number of
connected components which differ by a translation in $\RR^{r}$. One
of them is the tube $T:=F+i\RR^{r}$ with base
$$F=\big\{x\in\RR^{r}:2(\sin x_{1})^{2}+\sum_{k>1}e^{2x_{k}}
=1,\;|x_{1}|<\pi4\big\}\,.$$ With $2(\sin
x_{1})^{2}=1-\cos2x_{1}$ it is clear that $F$ is linearly equivalent
in $\RR^{r}$ to
$$\Pi_{+}:=\big\{x\in\RR^{r}:\cos
x_{1}=\sum_{k>1}e^{x_{k}},\,|x_{1}|<\pi/2\big\}$$ from \Lit{DAYA}. Here
$\PP_{r}\backslash O$ again is the union of $r{+}1$ projective
hyperplanes in general position, but all of them intersect $S$ and two
even tangentially. Figure 1 depicts in case $r=2$ the base of
$\Pi_{+}$ as the boundary of the `central' gray domain in $\RR^{2}$.
Also, the tube over the white region is the universal cover of
$S^{-}\cap O$, and the tube over every gray region is the universal
cover of $\{z\in B:z_{2}\ne0\}$ via $\phi$.

{\midinsert
 \epsfxsize=9true cm%
$$\hbox{\epsffile{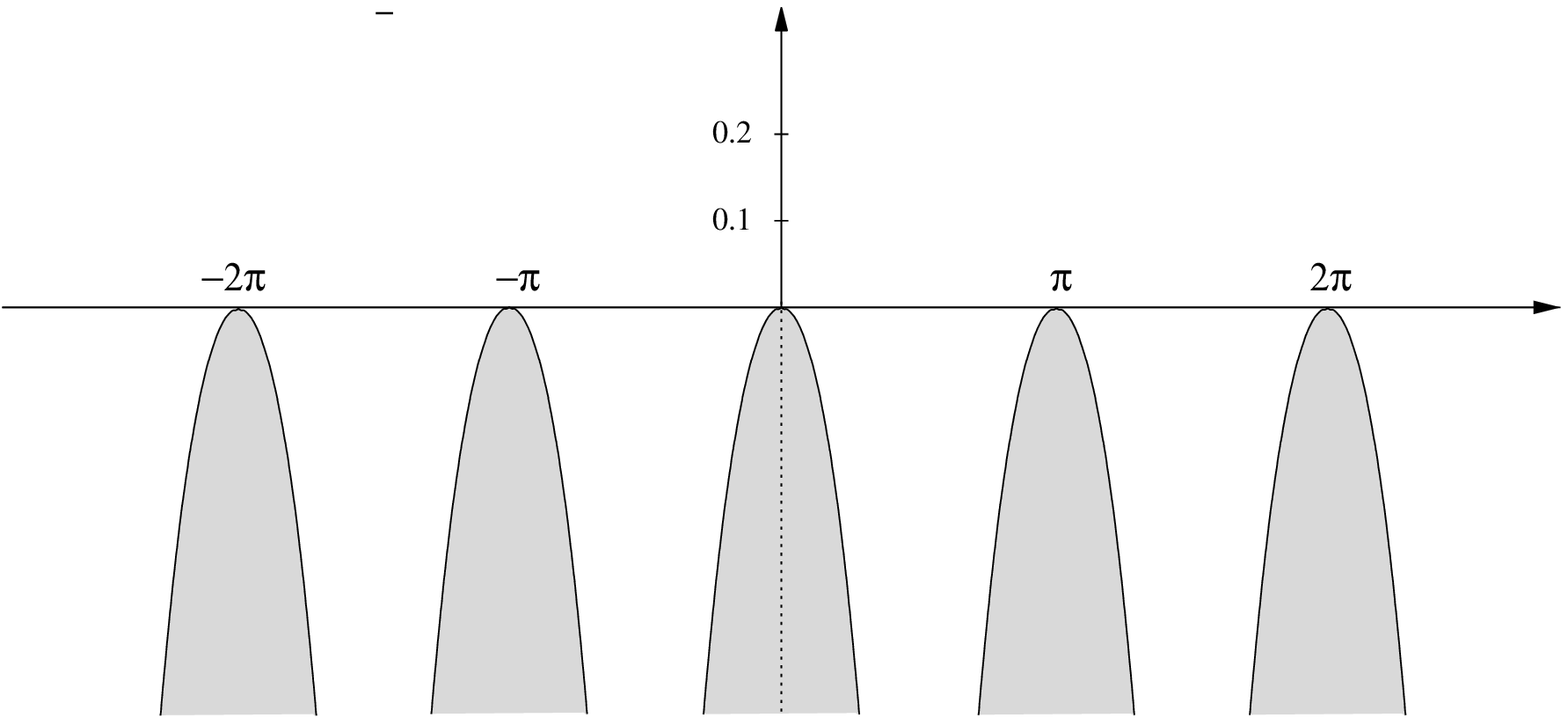}%
\Put{2}{-27}{$x_{1}$}%
\Put{-43.5}{-44}{$x_{2}$}
}$$
\vskip-5pt\centerline{\tt Figure 1}
\endinsert}

Notice that the abelian subalgebras $\7v\subset\7g$ giving the two
tube realizations $\Pi_{\pm}$ represent just the two conjugation
classes of {\sl Cartan subalgebras} of $\7g\cong\su(r,1)$ (= maximal
abelian subalgebras consisting of ad-semisimple elements).

\medskip To get further local tube realizations another description of
$S$ is convenient: Consider the classical Cayley transform
$\gamma\in\Aut(\PP_{r})$ defined by
$$\gamma([z]):=\big[z_{0}-z_{1},z_{0}+z_{1},
\sqrt2z_{2},\dots,\sqrt2z_{r}\big]\,.\Leqno{CA}$$ Then the
biholomorphic image $\gamma(S)$ in $\PP_{r}$ is of the form
$$S':=\gamma(S)=\big\{z\in\CC^{r}:z_{1}+\overline
z_{1}=\sum_{k=2}^{r}z_{k}\overline
z_{k}\big\}\;\cup\;\big\{[0,1,0,\dots,0]\big\}\,.$$ With $\7g=\hol(S)$
and $\7l=\7g^{\CC}$ as before let $\7g':=\hol(S')=\gamma_{*}\!\7g$.
For fixed $1\le s\le r$ let $\7v'$ be the linear span of the vector
fields\vskip3pt

\centerline{$i\dd{z_{1}}\!$, $\;iz_{r}\dd{z_{r}}$ and
$i(\dd{z_{j}}-z_{j}\dd{z_{1}}\!)\;\;$ for $\;1<r\le s$ and $s<j\le r$}
\noindent (written in the coordinate $z$ of $\CC^{r}$). Then $\7v'$ is
an abelian subalgebra of $\7g'$ and
$\7e':=\7v'\oplus\,i\!\7v'\subset\7l$ has the open orbit
$$O':=\{z\in\CC^{r}:z_{2}z_{3}\cdots z_{s}\ne0\}$$ in $\PP_{r}$. As a
consequence, $\PP_{r}\!\backslash O'$ is the union of $s$ mutually
different projective hyperplanes. As $\phi':\CC^{r}\to O'$ we can
choose
$$\phi'(z):=\Big((z_{1}-{1\over2}\sum_{j>s}z_{j}^{2}\,
,e^{z_{2}},\dots,e^{z_{s}},z_{s+1} ,\dots,z_{r})\Big)$$ and obtain the
corresponding tube realization with base
$$F_{s}:=\Big\{x\in\RR^{r}:x_{1}=\sum_{j=2}^{s}e^{2x_{j}}+
\sum_{j>s}x_{j}^{2}\Big\}\,.$$ $F_{s}$ is affinely equivalent to the
hypersurface $\Pi_{s-1,r-1}$ in \Lit{DAYA} and the tube $F_{s}+i\RR^{r}$
is the universal covering of
$$\{z\in S:(z_{1}-1)z_{2}\cdots z_{s}\ne0\}$$ via the map
$\phi:=\gamma^{-1}\phi'\,$. 

So far we have obtained $r+2$ local tube realizations of $S$ which are
mutually globally affinely inequivalent and closed in $\CC^{r}$. Among
these there is precisely one affinely homogeneous one -- the tube with
base $F_{1}=\{x\in\RR^{r}:x_{1}=\sum_{j>1}x_{j}^{2}\}$. This is the
unique algebraic tube realization and also the only case where
$\phi:\CC^{r}\to O$ is bijective and where $O\cap S$ is simply
connected.

By \Lit{DAYA} the examples above give, up to affine equivalence, all
closed smooth tube submanifolds in $\CC^{r}$ that are locally
CR-equivalent to the standard sphere $S=S_{1r}$.

\KAP{Further}{Further examples}

Our methods work best for CR-manifolds that are homogeneous (or at
least locally homogeneous). One way to get large classes of
CR-manifolds of this type is as follows: Choose a connected complex
Lie group $L$ acting holomorphically and transitively on a complex
manifold $Z$, that is, $Z=L/P$ for a closed complex Lie subgroup $P$
of $L$. Choose furthermore a real form $G$ of $L$, that is, a
connected real Lie subgroup such that $\7l=\7g^{\CC}$ for the
corresponding Lie algebras. Then for every $a\in Z$ the $G$-orbit
$S:=G(a)$ is an (immersed) CR-submanifold that is generically embedded
in $Z$ (since $\epsilon_{a}(\7l)=\T_{a}Z$). Clearly, the cases {\sl
$\;S$ open in $Z\;$} and {\sl $\;S$ totally real in $Z\;$} are not
interesting in our situation since for these the local CR-structure
is trivial and for every $a\in S$ there exists exactly one tube
realization of $(S,a)$ up to affine equivalence.

A case well understood in the group level is when $Z$ is a {\sl flag
manifold}, that is, $L$ is semisimple and $P$ is a parabolic
subgroup. Then, in particular, $Z$ is a compact rational projective
variety. The simplest flag manifold is the complex projective space
$\PP_{r}$ of dimension $r\ge1$. In this case we may take
$L=\SL(r{+}1,\CC)$ which is the universal cover of the group
$\Aut(\PP_{r})$. The only real forms $G$ of $L$ having an orbit in
$\PP_{r}$ that is neither open nor totally real are, up to
conjugation, the subgroups $\SU(p,q)$ with $p\ge q\ge1$ and
$m:=p+q=r+1$.  For the sake of completeness note that the real form
$G=\SL(m,\RR)$ has as unique non-open orbit the real projective space
$\PP_{r}(\RR)\subset Z$. This orbit is totally real and admits up to
affine equivalence a unique closed local tube realization in
$E=\CC^{r}$, namely $\RR^{r}\subset\CC^{n}$. The real form $\SU(m)$
and, in case $m$ is even, also the real form $\SL(m/2,\HH)$ act
transitively on $\PP_{r}$.

$\SU(1,1)$ is conjugate to $\SL(2,\RR)$ in $L$, so we assume $r>1$ in
the following. Then $G$ has again a unique non-open orbit in $Z$, the
compact hypersurface $S=S_{pq}$, compare \ruf{TU}. With
$\gamma\in\Aut(\PP_{r})$ the Cayley transform defined in \Ruf{CA}
$$Q:=\gamma(S)\cap\CC^{r}=\Big\{z\in\CC^{r}:z_{1}+\overline
z_{1}=\sum_{j=2}^{r}\epsilon_{j}z_{j}\overline
z_{j}\Big\},\quad\epsilon_{j}:=\cases{-1&$j\le p$\cr1&$j>
p$\cr}\,,\Leqno{WQ}$$ is the non-singular hyperquadric with Levi form
of type $(p-1,q-1)$. Now fix an integer $d$ with $1\le d\le r$. The
biholomorphic automorphism
$$(z_{1},\dots,z_{r})\longmapsto\Big(z_{1}+{1\over2}\sum_{j=2}^{d}
\epsilon_{j}z_{j}^{2},\,z_{2},\dots,z_{r}\Big)$$
of $\CC^{r}$ maps $Q$ to the submanifold
$$Q':=\Big\{z\in\CC^{r}:z_{1}+\overline
z_{1}={1\over2}\sum_{j=2}^{d}\epsilon_{j}(z_{j}+\overline
z_{j})^{2}+\sum_{j=d+1}^{r}\epsilon_{j}z_{j}\overline
z_{j}\Big\}\,.$$
Notice that $Q'$ has Siegel form, compare Section \ruf{Siegel},
$$Q':=\{(v,w)\in\CC^{d}\oplus\CC^{r-d}:(v+\overline v)-F(w,w)\in
C\}\,,\Leqno{SI}$$ where
$F(w,w):=\big(\sum^{r-d}_{j=1}\epsilon_{d+j}w_{j}\overline
w_{j},0,\dots,0\big)\in\RR^{d}$ and
$$C:=\Big\{x\in\RR^{d}:x_{1}=
{1\over2}\sum^{d}_{j=2}\epsilon_{j}x_{j}^{2}\Big\}\,.$$ In particular,
$Q'$ is a tube manifold in case $d=r$.

\medskip The next class of flag manifolds, to which our methods can
easily by applied, is given by the irreducible compact hermitian
symmetric spaces $Z$. Let $L$ be the universal covering of the
connected identity component of $\Aut(Z)$. Then $L$ is a simple
complex Lie group acting transitively on $Z$ and every real form of
$L$ has finitely many orbits in $Z$ that are all generically embedded
CR-submanifolds. There exists a real form $G$ of $L$ with an open
orbit $D$ that is biholomorphically equivalent to a bounded symmetric
domain. Suppose that $D$ is of tube type and choose a $G$-orbit
$S\subset Z$ that is neither open nor totally real. Then $S$ is Levi
degenerate (in fact is 2-nondegenerate) and $\hol(S)=\hol(S,a)$ is
isomorphic to the Lie algebra $\7g$ of $G$ for every $a\in S$, compare
\Lit{KAZT}. As a special example consider for fixed $p\ge2$ and
$m:=2p$ the Grassmannian $Z$ of all linear subspaces of dimension $p$
in $\CC^{m}$. Then $Z$ has complex dimension $n:=p^{2}$,
$L=\SL(m,\CC)$ and we can take $G=\SU(p,p)$. Now let $E:=\CC^{p\times
p}$ be the space of all complex $p{\times}p$-matrices and $V:=\{z\in
E:z^{*}=z\}$ the $\RR$-linear subspace of all hermitian matrices,
where $z^{*}$ is the transpose conjugate of the matrix $z$. The
$G$-orbits in $Z$ are in 1-1-correspondence to the cones
$$C^{p}_{j,k}:=\big\{x\in V:x\steil{has type}(j,k)\big\},\quad
j,k\ge0\steil{and}j+k\le p\,,\Leqno{OC}$$ in such a way that for every
$G$-orbit $S$ with corresponding $C^{p}_{j,k}$ the tube submanifold
$$T^{p}_{j,k}:=V+iC^{p}_{j,k}\;\subset\; E\Leqno{ST}$$ is
CR-equivalent to an open dense subset of $S$, see \Lit{KAZT}. Notice
that $T^{p}_{0,0}$ is the only closed tube submanifold of $E$ among
the $T^{p}_{j,k}$ in \Ruf{ST} and corresponds to the unique closed
$G$-orbit in $Z$ (totally real and diffeomorphic to the unitary group
$\U(p)$). On the other hand, all non-open tubes $T^{p}_{j,k}$, that is
$j+k<p$, are their own globalization in the sense of Section
\ruf{coarser}. Every cone $C^{p}_{j,k}$ is an orbit of the group
$\GL(n,\CC)$ acting linearly on $V$ by $x\mapsto gxg^{*}$, that is,
every tube $T^{p}_{j,k}$ is affinely homogeneous. All tubes
$T^{p}_{j,k}$ with $0<j+k<p$ satisfy the conditions {\rm P} and {\rm
Q} of Section \ruf{The}.

\KAP{Siegel}{CR-manifolds of Siegel type}

In the following we generalize the notion of a tube CR-manifold. Let
$V$ be a real and $W$ a complex vector space each of finite
dimension. Let furthermore $F:W\times W\to V^{\CC}$ be a $V$-hermitian
(vector valued) form, that is, complex linear in the first, antilinear
in the second variable and $F(w,w)\in V$ for every $w\in W$.
Throughout we assume that $F$ is nondegenerate, that is, $F(w,W)=0$
implies $w=0$ for every $w\in W$. For every real-analytic submanifold
$C\subset V$ and $\Im(x+iy):=y$ for all $x,y\in V$ then
$$\Sigma:=\{(z,w)\in V^{\CC}\oplus W:\Im z-F(w,w)\;\in\;C\}\Leqno{SG}$$ is
a real-analytic CR-submanifold of $E:=V^{\CC}\oplus W$ and is called a
{\sl Siegel CR-submanifold}. The CR-geometry of $\Sigma$ is closely related
to the {\sl associated tube} $T:=\Sigma\cap V^{\CC}=V+iC$ in $V^{\CC}$.
The submanifold $\Sigma$ is generically embedded in $E$ and $\Aut(\Sigma)$
contains the nilpotent subgroup
$$N:=\big\{\big(z,w\big)\mapsto
\big(z+v+2iF(w,c)+iF(c,c),w+c\big):v\in V,c\in W\big\}$$ acting by
affine transformations on $E$. Obviously $\Sigma=N(T)$ if we consider $T$
in the canonical way as submanifold of $\Sigma$. The Lie algebra
$$\7n=\{(2iF(w,c)+v)\dd z+c\dd w:v\in V,c\in W\}\subset\aut(\Sigma)$$ of
$N$ is nilpotent of step $\le2$ and can be considered as a subalgebra
of $\hol(\Sigma,a)$ with $\epsilon_{a}(\7n^{\CC})=E$ for every $a\in \Sigma$.

In a way, the nilpotent Lie subalgebras $\7n\subset\hol(M,a)$ play the
same role for Siegel realizations of a CR-manifold germ $(M,a)$ as the
abelian subalgebras $\7v\subset\hol(M,a)$ do for tube realizations.

\medskip Next we are interested in finite dimensionality conditions
for $\7g:=\hol(\Sigma,a)$, where $\Sigma$ is as in \Ruf{SG}. We start by
recalling (see e.g. \Lit{KAZT} for details) the\nline {\bf Iterated
Levi kernels.} Let $M$ be a CR-manifold of {\sl constant degeneracy}
(for instance if $M$ is locally homogeneous). Then there exists an
infinite descending chain of complex subbundles $$\H
M=\H^{\,0}\!M\supset\H^{\,1}\!M\supset\cdots\supset
\H^{\,k}\!M\supset\dots$$ where for every $a\in M$ the fiber
$\H_{a}^{\,k}\!M$, the $k$\th Levi kernel at $a$, is defined
recursively as follows: Choose a subset $\Xi\subset\Gamma(M,\H M)$
with $\epsilon_{a}(\Xi)=\H_{a}M$, where $\Gamma(M,\H M)$ is the space
of all smooth sections in $\H M$ over $M$. For every
$\eta\in\Gamma(M,\H^{\,k}\!M)$ the vector
$\eta_{a}\in\H_{a}^{\,k}\!M$ is in $\H_{a}^{\,k+1}\!M$ if and only if
$$[\xi,\eta]_{a}+i[\xi,i\eta]_{a}\,\in\,\H_{a}^{\,k}\!M\Steil{for
all}\xi\in\Xi$$ (this condition does not depend on the choice of
$\Xi$). In particular, $M$ is $k$-nondegenerate at every
point if and only if $\H^{\,k}\!M=0$, and $k$ is minimal with respect
to this property.

\Lemma{UD} Suppose that $\Sigma$ from \Ruf{SG} as well as the associated
tube $T=\Sigma\cap V^{\CC}$ have constant degeneracy. Then for every $a\in
T\subset \Sigma$ and every $k\ge0$ there exists a complex linear subspace
$W^{k}_{a}\subset W$ with $\H_{a}^{\,k}\!\Sigma\;=\;\H_{a}^{\,k}T\oplus
W^{k}_{a}$. Furthermore, $W^{0}_{a}=W$ and
$F(W^{k+1}_{a},W)\subset\H^{k}_{a}T$.  In particular, $H^{k}T=0$
implies $H^{k+1}\Sigma=0$.

\Proof We extend every $\xi\in\Gamma(T,\T \Sigma)$ to a smooth vector field
$\tilde\xi\in\Gamma(\Sigma,\T \Sigma)$ by requiring that for every $c\in W$ and
$\gamma\in N$ defined by $\gamma(z,w)=(z+2iF(w,c)+iF(c,c),w+c)$ we
have $\tilde\xi_{\gamma(z,0)} =d\gamma_{z}(\xi_{z})$ for all $z\in
T$. If we write $$\xi=f(z)\dd z+g(z)\dd w$$ with suitable smooth
functions $f:T\to V^{\CC}$ and $g:T\to W$, a simple computation
shows
$$\tilde\xi=(f(z-iF(w,w))+2iF(g(z),w))\dd z+g(z)\dd w\,.$$ From the
construction it is clear that $\xi\in\Gamma(T,H^{k}\Sigma)$ implies
$\tilde\xi\in\Gamma(\Sigma,H^{k}\Sigma)$ for all $k\in\NN$. Every
$\xi\in\Gamma(T,\T \Sigma)$ has a unique decomposition
$\xi=\xi^{1}+\xi^{2}$ with $\xi^{1}\in\Gamma(T,\T T)$ and
$\xi^{2}\in\Gamma(T,T\times W)$. Let $\Xi$ be the space of all
$\tilde\xi\in\Gamma(\Sigma,\H \Sigma)$ where $\xi\in\Gamma(T,\H \Sigma)$ has constant
second part $\xi^{2}$, that is, $\xi^{2}=c\dd w$ for some constant
vector $c\in W$. Then $\epsilon_{a}(\Xi)=H_{a}\Sigma$ is obvious. For $k=0$
the claim is obvious. Therefore assume as induction hypothesis that
the claim already holds for some fixed $k\ge0$.\nline Fix an arbitrary
$\eta\in\Gamma(T,\H^{k}\Sigma)$. Then $\eta_{a}=(\alpha,\beta)$ with
$\alpha\in\H^{k}_{a}(T)\subset V^{\CC}$ and $\beta\in W$. A simple
calculation shows that for every section $\xi=h(z)\dd z+c\dd
w\in\Gamma(T,\H T)$ with $\tilde\xi\in\Xi$ there exists a vector $e\in
W$ such that
$$[\tilde\xi,\tilde\eta]_{a}+i[\tilde\xi,i\tilde\eta]_{a}=
\big([\xi^{1},\eta^{1}]_{a}+i[\xi^{1},i\eta^{1}]_{a}-4iF(c,\beta)\dd
z\big)+e\dd w\,.\leqno{(*)}$$ Since $h(z)$ and $c$ can be chosen
independently for $\xi$ we derive from $(*)$ and the induction
hypothesis that $(\alpha,\beta)\in\H^{k+1}_{a}\Sigma$ implies
$\alpha\in\H^{k+1}_{a}T$ and $F(W,\beta)\subset\H^{k}_{a}T$.  Now
consider conversely an arbitrary $\alpha\in\H^{k+1}_{a}T$ and fix an
$\eta\in\Gamma(T,\H^{k+1}T)$ with $\eta_{a}=\alpha$. Then $(*)$ holds
with $\beta=e=0$ for every $\xi$ with $\tilde\xi\in\Xi$, that is,
$\alpha\in\H^{k+1}_{a}\Sigma$. This completes the induction step $k\to
k{+}1$.\qed

\smallskip As an application we state

\Proposition{HD} Let $\Sigma$ be an arbitrary Siegel submanifold as in
\Ruf{SG} and $T$ the associated tube manifold. Then \0 $\Sigma$ is
holomorphically nondegenerate if $T$ has the same property. \1 $\Sigma$ is
of finite type if $T$ has the same property or, if the set $F(W,W)$
spans the vector space $V^{\CC}$.

\Proof (i) Assume that $T$ is holomorphically nondegenerate. To show
that $\Sigma$ is holomorphically nondegenerate we only have to show that
$\Sigma$ is holomorphically nondegenerate at some point, compare Theorem
11.5.1 in \Lit{BERO}. We may therefore assume without loss of
generality that $T$ is of constant degeneracy and that
$\H^{k}T=0$. But then, as a consequence of Lemma \ruf{UD}, there
exists a domain $U\subset \Sigma$ of constant degeneracy with
$\H^{k+1}U=0$.\nline (ii) In a first step assume that $T$ is of finite
type in $a\in T$. Then the vector fields in $\Gamma(T,\T T)$ together
with all their iterated brackets span the tangent space $\T_{a}T$. For
all $\xi,\eta\in\Gamma(T,\T T)$ we have
$\tilde{[\xi,\eta]}=[\tilde\xi,\tilde\eta]$, where the extensions
$\tilde\xi,\tilde\eta\in\Gamma(\Sigma,\T \Sigma)$ are defined as in the proof of
\ruf{UD}. This shows that also $\Gamma(\Sigma,\H \Sigma)$ together with its
iterated brackets spans the tangent space $\T \Sigma_{a}$. From $N(T)=\Sigma$ we
get this property at every point of $\Sigma$.\nline Next assume that
$F(W,W)$ spans $V^{\CC}$. For every $c\in W$ and $\xi:=c\dd
w\in\Gamma(T,\H \Sigma)$ then
$$\tilde\xi,\tilde{i\xi}\in\Gamma(\Sigma,\H
\Sigma)\steil{and}\big[\tilde\xi,\tilde{i\xi}\big]=-4F(c,c)\dd z\,.$$
Since, by assumption, the vectors $F(c,c)$ span $V$, $\Sigma$ is of finite
type at every point of $T$ and hence also at every point of $\Sigma$.\qed

\KAP{Some}{Some Siegel CR-manifolds coming from bounded symmetric domains}

Irreducible bounded symmetric domains come in six types and for all
types the following considerations could be carried out in a uniform
(but more involved) approach. For simplicity we restrict our attention
only to the first type and there only to those domains that are not of
tube type: Fix integers $q>p\ge1$ and denote by $Z:=\2G_{p,q}$ the
Grassmannian of all $p$-dimensional linear subspaces in $\CC^{n}$,
$n:=p+q$. Then $Z$ is a compact complex manifold of complex dimension
$pq$, on which the complex Lie group $L:=\SL(n,\CC)$ acts transitively
by holomorphic transformations in a canonical way. Because of our
assumption $p\ne q$ the automorphism group $\Aut(Z)$ is connected and
has $L$ as universal cover. The real form $G:=\SU(p,q)$ of $L$ has
${p+2\choose2}$ orbits in $Z$. These can be indexed as
$M^{p,q}_{j,k}$, where $j,k\ge0$ are integers with $j+k\le p$. Indeed,
choose a $G$-invariant hermitian form $\Psi$ of type $(p,q)$ on
$\CC^{n}$ and let $M^{p,q}_{j,k}\subset Z$ be the set of all linear
subspaces, on which $\Psi$ has type $(j,k)$. For instance, the open
orbit $M^{p,q}_{p,0}$ is a bounded symmetric domain biholomorphically
equivalent to the operator ball
$$B:=\{z\in\CC^{p\times q}:(\One-zz^{*})\hbox{ positive
definit}\}\,,\Leqno{ZX}$$ where the matrix space $\CC^{p\times q}$ is
embedded in $Z$ as open dense subset by identifying every
$c\in\CC^{p\times q}$ with its graph
$\{(x,xc):x\in\CC^{p}\}\subset\CC^{n}$. In this way $M^{p,q}_{0,0}$,
the unique closed $G$-orbit in $Z$, corresponds to the extremal
boundary of $B$
$$\partial_{e}B:=\{z\in\CC^{p\times q}:\One=zz^{*}\}\,,$$ and
coincides also with the Shilov boundary of $B$. Notice that this
compact orbit already occurs as $S^{p,q}_{p}$ in Section
\ruf{Classification}.  Using a suitable Cayley transformation
$\gamma\in\Aut(Z)$ it can be shown that every $\gamma(M^{p,q}_{j,k})$
in the coordinate neighbourhood $\CC^{p\times q}\subset Z$ is the
CR-submanifold of Siegel type
$$\Sigma^{p,q}_{j,k}:=\{(z,w)\in\CC^{p\times
p}\oplus\CC^{p\times(q{-}p)}:\Im z-ww^{*}\in
C^{p}_{j,k}\}\,,\Leqno{ZY}$$ where $\Im z=\Quot{(z-z^{*})}{2i}$ and
the cone $C^{p}_{j,k}$ is as in \Ruf{OC}. For $V:=\{z\in\CC^{p\times
p}:z=z^{*}\}$ and $W:=\CC^{p\times(q{-}p)}$ the map $F:W\times W\to
V^{\CC}$, $(v,w)\mapsto vw^{*}$, satisfies $F(w,w)=0$ only for $w=0$
and its image $F(W,W)$ contains all rank-1-matrices in $\CC^{p\times
p}$. In particular, $F(W,W)$ spans $V^{\CC}$. Therefore, by
Proposition \ruf{HD}, all Siegel manifolds \Ruf{ZY} and hence all
$G$-orbits in $Z$ are of finite type. Now fix a $G$-orbit
$M=M^{p,q}_{j,k}\subset Z$ that is not open in $Z$, that is,
$j+k<p$. Denote by $T\subset\CC^{p\times p}$ the tube over
$C^{p}_{j,k}$. Then, if $j=k=0$ the tube $T$ is totally real and hence
$M\cong\partial_{e}B$ is Levi nondegenerate. In all other cases, that
is $0<j+k<p$, the tube $T$ is $2$-nondegenerate, compare Theorem 4.7
in \Lit{KAZT}. This implies with Lemma \ruf{UD} that every such $M$ is
Levi degenerate but is holomorphically nondegenerate. In particular,
for every non-open $G$-orbit $M$ in $Z$ and every $a\in M$ the Lie
algebra $\hol(M,a)$ has finite dimension and contains the simple Lie
algebra $\7g:=\su(p,q)$. On the other hand, since $G$ has a bounded
symmetric domain as orbit, for every $a\in M$ there is a local
coordinate $z$ around $a\in Z$ such $a$ is given by $z=0$ and that
$\7g^{\CC}$ contains all translation vector fields $c\dd z$ as well as
the Euler field $z\dd z$. With Proposition 3.1 in \Lit{KAZT} it
follows $\hol(M)=\hol(M,a)=\su(p,q)$ for every $a\in M$ and every
$G$-orbit $M$ in $Z$ which is neither open nor closed in $Z$.

\Proposition{PU} Every $G$-orbit $M\subset Z$ satisfies Condition {\rm
Q} of Section \ruf{The}. In case $M$ is neither open nor closed in $Z$
also Condition {\rm P} is satisfied.

\Proof The antilinear involution $z\mapsto \overline z$ of
$\CC^{p\times q}$ leaves the ball $B$ in \Ruf{ZX} invariant and
extends to an antiholomorphic involution $\tau$ of
$Z=\2G_{p,q}$. Therefore, $\tau$ leaves invariant every $G$-orbit in
$Z$. Now assume that the $G$-orbit $M$ is neither open nor closed in
$Z$. Then $\7g:=\hol(M)\cong\su(p,q)$ and for every $a\in M$ the
canonical restriction mapping $\rho_{a}:\7g\to\hol(M,a)$ is an
isomorphism of Lie algebras. For every $a\in M$ denote by
$\7g_{a}:=\{\xi\in\7g:\xi_{a}=0\}$ the isotropy subalgebra at $a$.  By
Proposition 2.11 in \Lit{KAPP}, $\7g_{a}=\7g_{b}$ for $a,b\in Z$ only
holds if $a=b$. The group $\Aut(M)\cong\PSU(p,q)$ is connected and for
$H:=\Aut(M)\cup\Aut(M)\tau$ the homomorphism $\Ad:H\to\Aut(\7g)$ is an
isomorphism, compare Proposition 4.5 in \Lit{KAPP}. In particular,
$\Aut(\7g)$ has two connected components. Now suppose that
$\phi:(M,a)\to(M,b)$ is either a CR-isomorphism or an
anti-CR-isomorphism of germs, where $a,b\in M$ are arbitrary
points. Then $\rho^{-1}_{b}\phi_{*}\rho_{a}$ is in $\Aut(\7g)$. In
case $\rho^{-1}_{b}\phi_{*}\rho_{a}$ is contained in the connected
identity component $\Int(\7g)$ of $\Aut(\7g)$ there exists a
transformation $g\in G$ such that $\rho^{-1}_{c}\psi_{*}\rho_{a}=\id$
for $c:=g(b)$ and $\psi:=g\,\phi:(M,a)\to(M,c)$. This implies $a=c$
and even $\psi=\id$ since $\rho^{-1}_{c}\psi_{*}\rho_{a}$ leaves
invariant all isotropy subalgebras $\7g_{x}$ for all $x\in M$ near
$a$. As a consequence, $\phi$ extends to the global transformation
$g^{-1}\in G$ in case $\rho^{-1}_{b}\phi_{*}\rho_{a}\in\Int(\7g)$.
But the case $\rho^{-1}_{b}\phi_{*}\rho_{a}\notin\Int(\7g)$ cannot
occur since otherwise $\rho^{-1}_{e}(\tau\phi)_{*}\rho_{a}\in\Int$ for
$e:=\tau(b)$ by the above reasoning would imply that $\tau\phi$ is a
CR-mapping, or equivalently, that $\phi$ is anti-CR.\qed

By the above considerations we know that for every non-open $G$-orbit
$M=M^{p,q}_{j,k}$ in $Z$ there is an integer $1\le k\le 3$ such that
$M$ is $k$-nondegenerate. In case $j+k=0$ we have $k=1$, and we claim
that $k=2$ in all other cases (compare also \Lit{FELS}): Indeed,
instead of $M$ we consider the Siegel manifold $\Sigma=\Sigma^{p,q}_{j,k}$ with
$V^{\CC}=\CC^{p\times p}$, $W=\CC^{p\times (q{-p})}$ and
$F(w,w)=ww^{*}$, compare \Ruf{ZY}. With $\rho:=j+k$ we write all $z\in
V^{\CC}$ and $w\in W$ as block matrices
$$z=\pmatrix{z_{11}&z_{12}\cr
z_{21}&z_{22}\cr}\steil{and}w=\pmatrix{w_{1}\cr w_{2}\cr}\;,$$ where
$z_{11}\in\CC^{\rho\times\rho}$, $w_{1}\in\CC^{\rho\times (q{-p})}$
and so forth.  Fix an element $a\in T=\Sigma\cap V^{\CC}$ with $a_{rs}=0$
for $(r,s)\ne(1,1)$. Then it is known that $$H^{k}_{a}T=\{z\in
V^{\CC}:z_{rs}=0\steil{if}k+r+s>3\}\,,$$ compare \Lit{KAZT} \p480.
This implies $w_{2}=0$ for every $w\in W^{1}_{a}$ and thus
$W^{2}_{a}=0$, that is, $H^{1}_{a}\Sigma\ne0$ and $H^{2}_{a}\Sigma=0$.

\medskip The antiholomorphic involution $\theta$ of $Z$ given on
$E=\CC^{p\times q}\subset\2G_{p,q}$ by $\theta(z)=-\overline z$ leaves
every Siegel manifold $\Sigma=\Sigma^{p,q}_{j,k}$ in \Ruf{ZY} invariant and has
fixed points there.  For every such fixed point $a\in \Sigma$ then
$\T_{a}^{-\theta}\Sigma=\RR^{p\times q}$, that is, \ruf{KR}.iii holds in this
situation. Assuming in the following that $\Sigma$ is not open in $E$ we
can have a local tube realization of $(\Sigma,a)$ associated with the
involution $\theta$ only if there is a maximal abelian subalgebra of
$\7g=\su(p,q)$ with dimension $pq$.  It can be shown that this is not
possible if $p>1$.

\vskip7mm {\gross\noindent References} \medskip
{\klein
\parindent 15pt\advance\parskip-1pt

\Ref{BART}Baouendi, M.S., Rothschild, L.P., Treves F.: CR structures with group action and extendibility of CR functions. Invent. math. {\bf 82}, 359-396 (1985).
\Ref{BAJT}Baouendi, M.S., Jacobowitz, H., Treves F.: On the analyticity of CR mappings. Ann. of Math. {\bf 122}, 365-400 (1985).
\Ref{BERO}Baouendi, M.S., Ebenfelt, P., Rothschild, L.P.: {\sl Real Submanifolds in Complex Spaces and Their Mappings}. Princeton Math. Series {\bf 47}, Princeton Univ. Press, 1998.
\Ref{BARZ}Baouendi, M.S., Rothschild, L., Zaitsev, D.: Equivalences of real submanifolds in complex space. J. Differential Geom. {\bf 59}, 301-351 (2001). 
\Ref{BOGG}Boggess, A.: {\sl CR manifolds and the tangential Cauchy Riemann complex.} Studies in Advanced Mathematics. Boca Raton, FL: CRC Press, 1991.
\Ref{CHMO}Chern, S.S., Moser, J.K.: Real hypersurfaces in complex manifolds. Acta. Math. {\bf 133}, 219-271 (1974).
\Ref{DAYA}Dadok, J.,Yang, P.: Automorphisms of tube domains and spherical hypersurfaces. Amer. J. Math. {\bf 107}, 999-1013 (1985).
\Ref{FELS}Fels, G.: Locally homogeneous finitely nondegenerate CR-manifolds. Math. Res. Lett. {\bf 14}, 693-922 (2007).
\Ref{FEKA}Fels, G., Kaup, W.: CR-manifolds of dimension 5: A Lie algebra approach. J. Reine Angew. Math, {\bf 604}, 47-71 (2007).
\Ref{FKAU}Fels, G., Kaup, W.: Classification of Levi degenerate homogeneous CR-manifolds of dimension 5. Acta Math. {\bf 201}, 1-82 (2008). 
\Ref{FLKA}Fels, G., Kaup, W.: Classification of commutative algebras and tube realizations of hyperquadrics. arXiv:0906.5549
\Ref{ISMI}Isaev, A.V., Mishchenko, M.A.: Classification of spherical tube hypersurfaces that have one minus in the Levi signature form. Math. USSR-Izv. {\bf 33}, 441-472 (1989).
\Ref{ISAE}Isaev, A.V.: Classification of spherical tube hypersurfaces that have two minuses in the Levi signature form. Math. Notes {\bf 46}, 517-523 (1989).
\Ref{ISAV}Isaev, A.V.: Global properties of spherical tube hypersurfaces. Indiana Univ. Math. {\bf 42}, 179-213 (1993).
\Ref{ISKA}Isaev, A.V., Kaup, W.:  Regularization of Local CR-Automorphisms of Real-Analytic CR-manifolds. arXiv:0906.3079
\Ref{KAPP}Kaup, W.: On the holomorphic structure of $G$-orbits in compact hermitian symmetric spaces. Math. Z. {\bf 249}, 797-816 (2005).
\Ref{KAPW}Kaup, W.: CR-quadrics with a symmetry property. arXiv:0907.4648
\Ref{KAZT}Kaup, W., Zaitsev, D.: On local CR-transformations of Levi degenerate group orbits in compact Hermitian symmetric spaces. J. Eur. Math. Soc. {\bf 8}, 465-490 (2006).
\Ref{PALA}{Palais, R.~S.:} A global formulation of the Lie theory of transformation groups. Mem. AMS 1957.
\Ref{TANA}Tanaka, N.: On the pseudo-conformal geometry of hypersurfaces of the space of $n$ complex variables. J. Math. Soc. Japan {\bf 14}, 397-429 (1962).
\Ref{ZAIT}Zaitsev, D.: On different notions of homogeneity for CR-manifolds. Asian J. Math. {\bf 11}, 331-340 (2007).
\bigskip
}

\smallskip\openup-4pt\parindent0pt
\hbox to 3cm{\hrulefill}
G. Fels\par e-mail: gfels@uni-tuebingen.de\par\vskip8pt
W. Kaup\par Mathematisches Institut, Universit\"at
T\"ubingen,\par Auf der Morgenstelle 10,\par  72076 T\"ubingen,
Germany\par e-mail: kaup@uni-tuebingen.de\par

\closeout\aux\bye